\documentclass[10pt]{amsart}

\usepackage{eucal, mathrsfs}
\usepackage{amsmath,amssymb,amsthm,amsfonts,amscd,amsopn}
\usepackage[OT2,OT1,T1]{fontenc}
\usepackage{newcent}
\usepackage{xspace}
\usepackage{epsfig,epic,eepic,latexsym,color}
\usepackage[all]{xy}
\usepackage{comment} 

\usepackage{latexsym}

\usepackage{fullpage}

\newcommand{\field}[1]{\mathbf #1}
\newcommand{\mf}[1]{\mathfrak #1}
\newcommand{\mc}[1]{\mathcal #1}
\newcommand{\ms}[1]{\mathscr #1}
\newcommand{\widebar}[1]{\overline{#1}}

\newcommand{\R}{\field R}
\renewcommand{\L}{\field L}
\newcommand{\C}{\field C}
\newcommand{\F}{\field F}
\newcommand{\Z}{\field Z}

\newcommand{\N}{\field N}

\newcommand{\simto}{\stackrel{\sim}{\to}}
\newcommand{\eps}{\varepsilon}
\renewcommand{\phi}{\varphi}

\newcommand{\ram}{\operatorname{Ram}}

\renewcommand{\hom}{\operatorname{Hom}}
\newcommand{\Hom}{\operatorname{Hom}}
\newcommand{\shom}{\ms H\!om}

\newcommand{\rshom}{\mathbf{R}\shom}

\newcommand{\send}{\ms E\!nd}

\newcommand{\spec}{\operatorname{Spec}}
\newcommand{\spf}{\operatorname{Spf}}

\renewcommand{\P}{\field P}

\DeclareMathOperator{\Pic}{Pic}

\DeclareMathOperator{\supp}{Supp}
\DeclareMathOperator{\pr}{pr}

\newcommand{\m}{\boldsymbol{\mu}}

\newcommand{\G}{\field G} 

\newcommand{\sep}{\text{\rm sep}}

\renewcommand{\H}{\operatorname{H}}

\newcommand{\GL}{\operatorname{GL}}
\newcommand{\PGL}{\operatorname{PGL}}
\DeclareMathOperator{\SL}{\operatorname{SL}}
\DeclareMathOperator{\ext}{\operatorname{Ext}}

\DeclareMathOperator{\Sym}{Sym}

\DeclareMathOperator{\per}{per}
\DeclareMathOperator{\ind}{ind}

\renewcommand{\(}{(\!\hspace{0.03em}(}
\renewcommand{\)}{)\!\hspace{0.03em})}
\DeclareMathOperator*{\tensor}{\otimes}

\renewcommand{\N}{\operatorname{N}}
\DeclareMathOperator{\Tr}{\operatorname{Tr}}
\DeclareMathOperator{\rk}{\operatorname{rk}}

\newcommand{\surj}{\twoheadrightarrow}
\newcommand{\inj}{\hookrightarrow}
\newcommand{\id}{\operatorname{id}}

\DeclareMathOperator{\coker}{\operatorname{coker}}

\DeclareMathOperator{\End}{\operatorname{End}}
\DeclareMathOperator{\aut}{\operatorname{Aut}}
\DeclareMathOperator{\isom}{\operatorname{Isom}}
\DeclareMathOperator{\M}{\operatorname{M}}
\DeclareMathOperator{\Br}{\operatorname{Br}}

\DeclareMathOperator{\V}{\bf V}

\newcommand{\Tw}{\mathbf{Tw}}
\DeclareMathOperator{\mTw}{Tw}
\newcommand{\Sh}{\mathbf{Sh}}
\DeclareMathOperator{\mSh}{Sh}
\DeclareMathOperator{\B}{\operatorname{\sf B\!}}

\newtheorem{lem}{Lemma}[subsection]

\makeatletter
\@addtoreset{lem}{section}
\@addtoreset{lem}{subsection}
\@addtoreset{lem}{subsubsection}
\@addtoreset{lem}{paragraph}
\makeatother

\renewcommand{\thelem}{\ifnum\value{subsubsection}>0{\thesubsubsection.\arabic{lem}}\else{\ifnum\value{subsection}>0{\thesubsection.\arabic{lem}}\else{\thesection.\arabic{lem}}\fi}\fi}

\newtheorem{thm}[lem]{Theorem}
\newtheorem*{theorem}{Theorem}

\newtheorem*{conjecture}{Conjecture}
\newtheorem{conj}[lem]{Conjecture}
\newtheorem*{standard conjecture}{Standard Conjecture}
\newtheorem{prop}[lem]{Proposition}

\newtheorem{cor}[lem]{Corollary}

\newtheorem*{substandard conjecture}{Standard Conjecture}

\theoremstyle{definition}
\newtheorem{defn}[lem]{Definition}

\newtheorem{para}[lem]{\indent}

\theoremstyle{remark}
\newtheorem{remark}[lem]{Remark}
\newtheorem{rem}[lem]{Remark}

\newtheorem{notn}[lem]{Notation}

\theoremstyle{definition}

\theoremstyle{remark}

\author{Max Lieblich}
\address{Fine Hall, Washington Road, Princeton, NJ}

\newcommand{\ov}{\overline}
\newcommand{\til}{\widetilde}

\newcommand{\bZ}{\mathbb{Z}}
\newcommand{\bF}{\mathbb{F}}

\newcommand{\bQ}{\mathbb{Q}}

\newcommand{\cD}{\mathcal{D}}

\newcommand{\cO}{\mathcal{O}}

\newcommand{\fP}{\mathfrak{P}}
\newcommand{\fB}{\mathfrak{B}}

\theoremstyle{remark}
\newtheorem*{case}{Case}

\title{Period and index in the Brauer group of an arithmetic surface}

\begin{document}
\bibliographystyle{plain}
\maketitle
\begin{abstract}
  In this paper we introduce two new ways to split ramification of
  Brauer classes on surfaces using stacks.  Each splitting method
  gives rise to a new moduli space of twisted stacky vector bundles.
  By studying the structure of these spaces we prove new results on the
  standard period-index conjecture.  The first yields new bounds on
  the period-index relation for classes on curves over higher local
  fields, while the second can be used to relate the Hasse principle
  for forms of moduli spaces of stable vector bundles on pointed curves over
  global fields to the period-index problem for Brauer groups of
  arithmetic surfaces.  We include an appendix by Daniel Krashen
  showing that the local period-index bounds are sharp.
\end{abstract}

\tableofcontents

\section {Introduction}

Let $K$ be a field and $D$ a central division algebra over $K$.  There
are two basic numerical invariants one can attach to $D$.  First, it
is a standard fact that the $K$-dimension of $D$ is a square, say
$\dim_K D=n^2$, from which one can extract $n$, the \emph{index\/} of
$D$.  Second, $D$ gives rise to an element of the Brauer group of $K$,
which is a torsion group.  We can thus read off the order of $D$ in
$\Br(K)$, known as the \emph{period\/} of $D$.  Unless noted
otherwise, \emph{throughout this paper, we will only work with Brauer
  classes $\alpha\in\Br(K)$ of period prime to the characteristic
  exponent of $K$\/}.

The principal results of this paper are the following.  Given a
positive integer $M$, let $M'$ denote the submonoid of $\N$ consisting
of natural numbers relatively prime to $M$.
\begin{defn} 
  Let $M$ and $a$ be positive integers.  A field $k$ has
  \emph{$M'$-Brauer dimension at most $a$\/} if for every finitely
  generated field extension $L/k$ of transcendence degree $t$ at most
  $1$ and every Brauer class $\alpha\in\Br(L)$ with period contained
  in $M'$ one has $\ind(\alpha)|\per(\alpha)^{a-1+t}$.
\end{defn}
The case of particular interest is when $M=p$ for some prime $p$
(usually the characteristic exponent of $K$ or of a residue field of
$K$ under some valuation).
\begin{notn}
  Throughout this paper, when we use the phrase ``Brauer dimension at
  most $a$'' we will always tacitly assume that $a$ is a positive
  integer.
\end{notn}

\begin{theorem}[Theorem \ref{T:main}]
  Suppose $K$ is the fraction field of an excellent Henselian discrete
  valuation ring with residue field $k$ of characteristic exponent
  $p$.  If $k$ has $M'$-Brauer dimension at most $a$ then $K$ has
  $(Mp)'$-Brauer dimension at most $a+1$.
\end{theorem}

We can get an inductive corollary as follows.  By a \emph{$d$-local
  field\/} with residue field $k$ we will mean a sequence of fields
$k_d,\ldots,k_0=k$ such that for $i=1,\ldots,d$, $k_i$ is the fraction
field of an excellent Henselian discrete valuation ring with residue
field $k_{i-1}$.  Natural examples are given by iterated Laurent
series fields $k\( x_1\) \cdots\(x_d\)$.  

\begin{cor} If $K$ is a $d$-local field with residue field $k$ of
  $M'$-Brauer dimension at most $a$ and characteristic exponent $p$,
  then $K$ has $(Mp)'$-Brauer dimension at most $a+d$.
\end{cor}

A concrete application of this result yields various instances of the
standard conjecture (recalled below).

\begin{cor}\label{appendix-cor} 
  Let $L/K$ be a field extension of transcendence degree $1$ and
  $\alpha\in\Br(L)$ a Brauer class.
  \begin{enumerate}
  \item If $K=k\(x_1\)\cdots\(x_{d}\)$ is the field of iterated
    Laurent series over an algebraically closed field and
    $\per(\alpha)$ is relatively prime to the characteristic exponent
    of $k$ then $\ind(\alpha)|\per(\alpha)^{d}$.
  \item If $K$ is $d$-local with finite residue field of
    characteristic $p$ and $\per(\alpha)$ is realtively prime to $p$
    then $\ind(\alpha)|\per(\alpha)^{d+1}$.
  \item If $K$ is an algebraic extension of the maximal unramified
    extension of a local field of residue characteristic $p$ and
    $\per(\alpha)$ is relatively prime to $p$, then
    $\ind(\alpha)=\per(\alpha)$.
  \end{enumerate}
\end{cor}

Linking local global, we also prove the following. Let $K$ be a global
field and $C/K$ a proper, smooth, geometrically connected curve with a
rational point $P$.  For the purposes of this paper, we will refer to
the following statement as the \emph{Colliot-Th\'el\`ene
  conjecture\/}.  (It is a special case of a more general conjecture
which may be found in \cite{MR1743234}.)

\begin{conjecture}[Colliot-Th\'el\`ene]
  If $X$ is a smooth projective geometrically rational variety over a
  global field $K$ with geometric Picard group $\Z$, then $X$ has a
  $0$-cycle of degree $1$ if and only $X\tensor_K K_{\nu}$ has a $0$-cycle
  of degree $1$ for every place $\nu$ of $K$.
\end{conjecture}

\begin{theorem}[Theorem \ref{sec:curves-over-global-6}]
  Assuming the Colliot-Th\'el\`ene conjecture, any
  $\alpha\in\Br(K(C))$ of odd period such that 
  \begin{enumerate}
  \item $\alpha$ is unramified over primes dividing $\per(\alpha)$
  \item $\alpha|_P\in\Br(K)$ is trivial
  \end{enumerate}
satisfies
  $\ind(\alpha)|\per(\alpha)^2$.
\end{theorem}

When these surfaces are equicharacteristic, i.e., they are surfaces
over finite fields, then this is part of the standard conjecture for
$C_3$-fields.  Moreover, as shown in Corollary
\ref{sec:curves-over-global-4} below, one can deduce the standard
conjecture for classes of odd period over fields of transcendence
degree $2$ over finite fields from the above theorem.

The key idea underlying these results is the observation that one can
split ramification of Brauer classes using stacks and produce moduli
problems -- twisted sheaves over various kinds of stacky curves --
whose solution gives bounds on the period-index relation.  

After some stack-theoretic preliminaries in Section
\ref{S:stacky-mama}, we describe two methods for splitting
ramification by stacks in Section \ref{S:splitting on curves}.  We
study the Brauer groups of certain finite classifying stacks in
Section \ref{S:coho-mama} in preparation for a discussion of stacky
moduli problems in Section \ref{sec:residual}.  Finally, in Sections
\ref{S:local-mama} and \ref{sec:curves-over-global-1} we use the two
different stacky splitting methods to study the local and global
period-index results described above.  At the end is an appendix by
Daniel Krashen showing that our local results are sharp.

To put this work in context, let us make a few historical remarks and
remind the reader of standard conjectures on the relation between
period and index.  Basic Galois cohomology shows that
$\per(D)|\ind(D)$ and that both have the same prime factors, so that
$\ind(D)|\per(D)^{\ell}$ for some power $\ell$.  The exponent $\ell$
is still poorly understood, but the following ``conjecture'' has
gradually evolved from its study.  (The reader is referred to page 12
of \cite{ctl} for what appears to be the first recorded statement of
the conjecture, in the case of function fields of complex algebraic
varieties.)

\begin{standard conjecture}
  If $K$ is a ``$d$-dimensional'' field and $\alpha\in\Br(K)$, then
  $\ind(\alpha)|\per(\alpha)^{d-1}$.
\end{standard conjecture}

The phrase ``$d$-dimensional'' is left purposely vague.  It is perhaps
best to illustrate what it could mean by summarizing which cases of
the standard conjecture have been proven.  In \cite{Tsen}, Tsen proved 
that the Brauer group of the function field of a curve over an
algebraically closed field is trivial.  (More generally, his proof
shows that any $C_1$-field has trivial Brauer group.)  The next progress in this
direction was not made until de Jong proved in \cite{dejong-per-ind}
that the period and index are equal for Brauer classes over the
function field of a surface over an algebraically closed field.  (By
contrast, it is still unknown whether this equality holds for any
curve over a $C_1$-field, or for any $C_2$-field.)  

Threefolds have so far been quite difficult to handle, and in fact no one has
even proven that there is a bound on the period-index relation for
classes in $\Br(\C(x,y,z))$ (or the function field of any other
particular threefold).  The methods used to prove
these theorems have no inductive structure as one increases the
dimension of the ambient variety; this is in part responsible for the
difficulties faced in making progress along these lines.
For example, de Jong's proof makes essential use of the deformation theory of
Azumaya algebras on surfaces, something which one can prove becomes
far more complicated on an ambient threefold.

In a slightly different direction (and between Tsen and de Jong in
history), Saltman proved in \cite{saltman-fix} and \cite{saltman} that
$\ind | \per^2$ for curves over $p$-adic fields, and he provided a
method for showing that if $K$ is the fraction field of an excellent
Henselian discrete valuation ring with residue field $k$, and if $\ind
| \per^a$ for every Brauer class on every curve over $k$, then one has
that $\ind | \per^{a+2}$ for every Brauer class on every curve over
$K$.  His techniques provide for a kind of induction as one proceeds
up a ladder of Henselian fields, with the exponent in the relation
increasing by $2$ each time.

The key to Saltman's method is splitting the ramification of Brauer classes
by making branched covers, which allows one to reduce questions on
curves over the fraction field to curves over the residue
field. Because of the singularities which appear on branched covers,
in order to completely split the ramification of a Brauer class
Saltman is forced to make two covers, which is responsible for the
increase of the exponent by $2$.

The methods used in this paper are in some sense a hybrid of the ideas
of Saltman and de Jong: we make some branching (either with schemes or
stacks) and then use moduli theory to reduce a period-index relation
to a rational point problem.

Connections between the Hasse principle and the period-index problem
are also discussed in \cite{period-index-paper}, where this connection
is used to \emph{disprove\/} the Hasse principle (and the
Colliot-Th\'el\`ene conjecture) for smooth projective geometrically
rational varieties over $\C(x,y)$ with geometric Picard group $\Z$.
It should be interesting to investigate this connection in more detail
over global fields.  Some preliminary results in this direction are
described in \cite{lieblich-nutz}.

\section*{Notation}

We adopt the following notations:
\begin{enumerate}
\item Given a local ring $A$ with residue field $\kappa$ and a
  separable closure $\kappa\subset\kappa^\sep$, we let $A^{hs}$
  denote the strict Henselization of $A$ with respect to
  $\kappa^\sep$.  We may occasionally abuse the notation and leave
  the choice of $\widebar\kappa$ ambiguous.
\item Given an abelian group $G$, we let $G'$ denote the subgroup
  consisting of elements of order prime to $p$, where $p$ is the
  characteristic exponent of the base field.
\item Given a (possibly generic) point $p$ of a scheme $X$, we will
  write $\kappa(p)$ for the residue field of $p$.  If $X$ is integral,
  we will also write $\kappa(X)$ for the function field (the residue
  field at the generic point).
\end{enumerate}

Some additional stack-specific notation will be introduced in Section
\ref{S:stack-nots}.

\section*{Acknowledgments}

I would like to thank Jean-Louis Colliot-Th\'el\`ene, A.\ J.\ de Jong,
David Harbater, Daniel Krashen, David Saltman, and Fabrice Orgogozo
for numerous helpful discussions, observations, and corrections.

\section {Stack-theoretic preliminaries}
\label{S:stacky-mama}
\subsection{Notation and basic results}\label{S:stack-nots}

The following notations will be essential.

\begin{enumerate}
\item Given an algebraic space $X$, an invertible sheaf $\ms L$, and
  an integer $n$, the $\m_n$-gerbe of $n$th roots of $\ms L$ will be
  denoted $[\ms L]^{1/n}$.  This is a stack on $X$ which solves the
  moduli problem given by pairs $(\ms M,\phi)$ with $\ms M$ an
  invertible sheaf and $\phi:\ms M^{\tensor n}\simto\ms L$ is an
  isomorphism.  It is easy to check that the class of this gerbe in
  $\H^2(X_{\text{\rm fppf}},\m_n)$ is the image of the isomorphism
  class of $\ms L$ under the map $\Pic(X)\to\H^2(X_{\text{\rm
      fppf}},\m_n)$ induced by the Kummer sequence (in the fppf
  topology).  In particular, there is a $[\ms L]^{1/n}$-twisted
  invertible sheaf (see Section \ref{S:tw-sh} below)$\ms M$ with $\ms
  M^{\tensor n}\cong\ms L$.
\item A \emph{surface\/} is a quasi-compact excellent algebraic space of
  equidimension $2$.
\item An \emph{orbisurface\/} is a quasi-compact excellent 
  Deligne-Mumford stack of equidimension $2$ with trivial generic stabilizer.
\item Given an integer $n$, an \emph{$A_{n-1}$-orbisurface\/} (resp.\
  \emph{Zariski $A_{n-1}$-orbisurface\/}) is a separated regular tame
  orbisurface $\ms X$ with a coarse moduli space $\pi:\ms X\to X$ such
  that for each point $x\in X$, the localization $\ms
  X\times_X\spec\ms O_{x,X}^{hs}$ (resp.\ $\ms X\times_X\spec\ms
  O_{x,X}$) has the form $[R/\m_n]$, where $R$ is a regular local ring
  of dimension $2$ and $\m_n$ acts on $\mf m_R/\mf m_R^2$ by
  $\sigma\oplus\sigma^{\vee}$, where $\sigma$ is the standard
  character of $\m_n$.
\item Given a stack $\ms X$, the \emph{inertia stack\/} is the
  relative group-stack $\ms I(\ms X)\to\ms X$ representing the functor
  which to an object $\phi:X\to \ms X$ assigns the group $\aut(\phi)$
  of automorphisms of $\phi$ in the fiber category of $\ms X$ over
  $X$.  If $\ms X$ is an Artin stack locally of finite presentation
  over an algebraic space $S$, then it is a standard fact that $\ms
  I(\ms X)\to\ms X$ is of finite presentation.
\item Given a stack $\ms X$, the \emph{stacky locus\/} is the
  stack-theoretic support of the inertia stack.
\item Given a Deligne-Mumford stack $\ms X$ with coarse moduli space
  $\pi:\ms X\to X$ and a point $p\in X$, the \emph{residual gerbe\/}
  at $p$ is the reduced structure on the preimage of $p$ in $\ms X$.
  (This is a special case of a more general notion described in
  section 11 of \cite{l-mb}.)  If $\ms X$ is an $A_{n-1}$-orbisurface,
  then the residual gerbe at any point $p$ is either isomorphic to
  $\spec\kappa(p)$ (which we will call \emph{trivial\/}) or is a
  $\m_n$-gerbe over $\kappa(p)$; for a Zariski $A_{n-1}$-orbisurface,
  the non-trivial residual gerbes have trivial cohomology class in
  $\H^2(\kappa(p),\m_n)$.
\item Given a residual gerbe $\xi$ lying over $p$ with residue field
  $\kappa$, we will write $\widebar{\xi}$ for the base change of $\xi$
  to an algebraic closure $\widebar{\kappa}$.  (As is often the case,
  we may abuse notation -- as we have already -- and leave the choice
  of embedding $\kappa\inj\widebar{\kappa}$ out of the notation.)
\item Given a Deligne-Mumford stack $\ms X$ with coarse moduli space
  $\pi:\ms X\to X$ and a residual gerbe $\xi\in\ms X$ over $p\in
  X$, we will call the base change $\ms X\times_X\spec\ms
  O^{h}_{p,X}$ (resp.\ $\ms X\times_X\spec\ms
  O^{hs}_{p,X}$) the \emph{Henselization\/} (resp.\ \emph{strict
    Hen\-selization\/}) of $\ms X$ at $\xi$.
\end{enumerate}
For basic results on gerbes, the reader is referred to \cite{giraud}
or, for a condensed discussion, \cite{period-index-paper} or
\cite{twisted-moduli}.

When blowing up $A_{n-1}$-orbisurfaces, we will need the following
generalization of their local structure:

\begin{defn}
  Given a regular local ring $A$ of dimension $2$, an action of
  $\m_{n}$ on $A$ \emph{has type $(a,b)\in\Z/n\Z^{\oplus 2}$\/} if the
  induced representation of $\m_{n}$ on $\mf m_A/\mf m_A^2$ splits
  as $\sigma^{\tensor a}\oplus\sigma^{\tensor b}$.  We will also say
  that the quotient stack $[\spec A/\m_{n}]$ has type $(a,b)$.
\end{defn}

Given any tame regular Deligne-Mumford surface, a residual gerbe
of the form $\B{\m_{n}}$ has a corresponding type: it is the type
of the Henselization of the stack at that gerbe.

\subsection {Twisted sheaves on a topos and deformation theory}
\label{S:tw-sh}

Let $(X, \ms O)$ be a ringed topos.  For the sake of convenience, we
will assume that there is a covering $U$ of the final object such that
$(U, \ms O _ {U})$ and $(U \times U, \ms O _ {U \times U})$ are
isomorphic to big fppf topoi of schemes.  (Among other things, this
hypothesis ensures that quasi-coherent sheaves are an abelian
category, and it is satisfied by all algebraic stacks.)  The
references for the results in this section are
\cite{period-index-paper} and \cite{twisted-moduli} and the references
therein.

Let $A$ be either the sheaf $\G_m$ of units or the sheaf $\m_n$ of
$n$th roots of unity, with $n$ assumed invertible in $\Gamma (X, \ms O
_ {X})$.  Let $\ms X$ be an $A$-gerbe on $X$, and let $\ms U$ be the
restriction to $U$.  For any sheaf $\ms F$ on $\ms X$, pullback by
automorphisms (i.e., sections of the inertia stack) defines a right
action $\ms F \times A \to \ms F$.

\begin{defn} An \emph{$\ms X$-twisted sheaf\/} is an $\ms O _ {\ms
    X}$-module $\ms F$ such that the natural action $\ms F \times A
  \to \ms F$ agrees with the right $A$-action induced by the (left) module
  structure on $\ms F$.
\end{defn}

As usual, one can define coherent and quasi-coherent twisted sheaves,
and one can speak of locally free (quasi-coherent) twisted sheaves.
The basic fact underlying the usefulness of twisted sheaves for the
period-index problem is the following, which may be found in
Proposition 3.1.2.1 of \cite{period-index-paper}.

\begin{lem}
  If $\ms X\to X$ is a $\G_m$-gerbe over a regular algebraic
  space of dimension at most $2$ then there is a locally free $\ms
  X$-twisted sheaf of rank $n$ if and only if the index of the class
  $[\ms X]\in\H^2(X,\G_m)$ divides $n$.  
\end{lem}

More concretely, Azumaya algebras on $X$ with Brauer class $[\ms X]$
are precisely those algebras which may be written as the direct images
of trivial algebras $\send(\ms V)$ on $\ms X$ with $\ms V$ a locally
free $\ms X$-twisted sheaf.

Thus, if we wish to show that the period equals the index, we need
only look for a locally free twisted sheaf whose rank equals the period;
this places a geometric overlay on the seemingly algebraic
period-index problem.  This apparently trivial fact has many useful
consequences, some of which we will see below (and various others of
which are described in \cite{period-index-paper}).

Let $0 \to
I \to \ms O \to \ms O _ 0 \to 0$ be a small extension of rings in $X$.
Given an $A$-gerbe $\ms X\to (X,\ms O)$, define the \emph{reduction\/}
of $\ms X$ to be the $A(\ms O_0)$-gerbe $\ms X_0\to X_0:=(X,\ms O_0)$
given by applying Proposition 3.1.5 of \cite{giraud} to the map $A(\ms
O)\to A(\ms O_0)$ of sheaves of abelian groups.  (This is just a
stacky rigidification of the pullback map on cohomology $\H^2(X,A(\ms
O))\to\H^2(X,A(\ms O_0))$ induced by the map of
ringed topoi $(X,\ms O_0)\to(X,\ms O)$.)

\begin{prop}\label{P:deformation theory}  With the above notation, 
  let $F _ 0$ be an $\ms X _ 0$-twisted sheaf.
  \begin{enumerate}
  \item There is an obstruction in $\ext_{\ms X_{0}} ^ 2 (F _ 0, I
    \tensor F _ 0)$ to deforming $F _ 0$ to an $I$-flat
    $X$-twisted sheaf.
  \item If $X$ is a scheme of dimension 1 and $F_{0}$ is locally
    free, then for any $I$ we have $$\ext ^ 2 (F _ 0, I \tensor F _ 0) = \H ^
    2 (X, \shom(F _ 0, I \tensor F _ 0)) = 0.$$
  \item If $X$ is an algebraic stack, the Grothendieck existence
    theorem holds for formal deformations of coherent twisted sheaves.
  \end{enumerate}
\end{prop}

For further details and proofs, we refer the reader to
\cite{twisted-moduli} and \cite{period-index-paper}.  The last
statement is not found any of the references; it follows from the
Grothendieck existence theorem for stacks \cite{olsson-starr} applied
to the results of Section \ref{S:stack-on-a-stack} below, which shows that
gerbes on big sites of algebraic stacks are representable by algebraic
stacks.

\subsection{Moduli spaces of twisted sheaves}
\label{sec:moduli}

We briefly summarize a few facts about the moduli space of twisted
sheaves on a curve.  The purpose of this section is mostly to fix
notation.  We refer the reader to \cite{twisted-moduli} for further
details and proofs.

Let $X/K$ be a smooth proper curve over a field and $\ms G\to X$ a
$\m_n$-gerbe.  Given $L\in\Pic(X)$, there is a moduli stack of
\emph{stable twisted sheaves\/} of rank $n$ and determinant $L$,
written $\Tw^s_{\ms G/K}(n,L)$.  It is a $\m_n$-gerbe over its coarse
moduli space $\mTw^s_{\ms G/K}(n,L)$.  Moreover, $\mTw^s_{\ms
  G/K}(n,L)$ is a form of the space of stable sheaves
$\mSh^s_{X/K}(n,L')$ for an appropriate $L'\in\Pic(X)$.  (In
particular, it arises as a twist by automorphisms in $\Pic_{X/K}[n]$.)
Thus, $\mTw^s_{\ms G/K}(n,L)$ is geometrically (separably)
unirational; if $\deg L$ is appropriately chosen, then it is in fact
proper and geometrically rational, with (geometric) Picard group $\Z$.

\subsection {A stack on a stack is a stack}\label{S:stack-on-a-stack}

We assume for
the sake of simplicity that all categories in this section admit
arbitrary (non-empty) coproducts and all finite fibered products
(whose indexing sets are non-empty), but not necessarily that they
have final objects.  Given a category $S$, we will interchangeably use
the terminology ``$S$-groupoid'' and ``category fibered in groupoids
over $S$.''

\begin{para}\label{Para:setup-stack}
Let $S$ be a category and $\pi_X:X\to S$ and $\pi_Y:Y\to S$ two
categories fibered in groupoids.
\end{para}
\begin{defn}\label{D:fib-rep}
  Given a functor $F:Y\to X$ of fibered categories, the \emph{fibered
    replacement\/} of $F$ is the fibered category $(Y,F)\to X$
  consisting of triples $(y,x,\phi)$ with $y\in Y$, $x\in X$, and
  $\phi:F(y)\simto x$ an isomorphism.  The arrows
  $(y,x,\phi)\to(y',x',\phi')$ are arrows $y\to y'$ and $x\to x'$ such that the
  diagram
$$\xymatrix{F(y)\ar[r]\ar[d] & x\ar[d]\\
F(y')\ar[r] & x'}$$ commutes.
\end{defn}

There is a natural functor $Y\to (Y,F)$ given by sending $y$ to
$(y,F(y),\id)$ and a natural functor $(Y,F)\to Y$ sending $(y,x,\phi)$
to $y$.

\begin{lem}\label{L:fib-lem} Given $F:Y\to X$ as above,
  \begin{enumerate}
  \item the natural functor $(Y,F)\to X$ is a category fibered in groupoids;
  \item the natural functor $(Y,F)\to S$ is a category fibered in groupoids;
  \item the natural functors $Y\to (Y,F)$ and $(Y,F)\to Y$ define an
    equivalence of fibered categories over $S$.
  \end{enumerate}
\end{lem}
\begin{proof}
  To check that $(Y,F)\to X$ is a category fibered in groupoids, suppose
  $(y,x,\phi)$ is an object of $(Y,F)$ and $\alpha:x'\to x$ is an
  arrow in $X$.  Since $Y$ is a fibered category, there is an arrow
  $\beta:y'\to y$ such that $\pi_Y(\beta)=\pi_X(\alpha)$.  Since every
  arrow in a category fibered in groupoids is Cartesian, there is an
  isomorphism $\phi':F(y')\to x'$ such that the diagram
$$\xymatrix{F(y')\ar[r]\ar[d] & x'\ar[d]\\
F(y)\ar[r] & x}$$
commutes.  But this says precisely that $(y',x',\phi')\to(y,x,\phi)$
lies over $\alpha$.  We will be finished once we check that every
arrow $(y',x',\phi')\to(y,x,\phi)$ is Cartesian.  Thus, suppose
$(y'',x',\phi'')\to(y,x,\phi)$ is another arrow such that the
underlying arrow $x'\to x$ is the same.  In this case, we know that
$y''\to y$ and $y'\to y$ have the same image in $S$, whence there is
a unique arrow $y''\to y'$ over $y$.  This results in a diagram
$$\xymatrix{&F(y'')\ar[dl]\ar[dr] & \\
F(y')\ar[d]\ar[rr] & & F(y)\ar[d]\\
x'\ar[rr] & & x.}$$
Since the arrow $F(y'')\to x'$ is uniquely determined by the arrows
$F(y'')\to F(y)$ and the arrow $x'\to x$, we see that there results a
map $(y'',x',\phi'')\to(y',x',\phi')$, as desired.

That $(Y,F)\to S$ is a category fibered in groupoids is a formal
consequence of the first part and the fact that $X\to S$ is a category
fibered in groupoids.

It is immediate that the composition $Y\to (Y,F)\to Y$ is equal (!) to
the identity.  In the other direction, the composition $(Y,F)\to Y\to
(Y,F)$ sends $(y,x,\phi)$ to $(y,F(y),\id)$.  The diagram
$$\xymatrix{F(y)\ar[r]^\phi\ar[d]_{\id} & x\ar[d]^{\phi^{-1}}\\
  F(y)\ar[r]_{\id} & F(y) }$$ shows that the map $(\id_y,\phi^{-1})$
is an isomorphism $(y,x,\phi)\simto(y,F(y),\id)$, which shows that the
composition of the two functors is naturally isomorphic to
$\id_{(Y,F)}$.  This gives the last statement.
\end{proof}

\begin{prop}\label{P:fib-equiv}
  Given a category fibered in groupoids $X\to S$, the forgetful
  functor and the functor $\nu:(F:Y\to
  X)\mapsto(Y,F)$ give an equivalence between the slice 2-category of
  $S$-groupoids over $X$ and the 2-category of $X$-groupoids.
\end{prop}
\begin{proof}
  It is clear by construction that the (strong) equivalence
  $Y\cong(Y,F)$ of $S$-groupoids proven in Lemma \ref{L:fib-lem} is
  functorial in $Y$ and $F$.  In particular, it follows that given two
  functors $F:Y\to X$ and $F':Y'\to X$, the functor
  $\hom_{X}(Y,Y')\to\hom_X((Y,F),(Y',F'))$ is an equivalence of
  groupoids, so that $\nu$ is fully faithful (in the strong sense).
  On the other hand, any $X$-groupoid $H:Z\to X$ is naturally an
  $S$-groupoid (as $X$ is an $S$-groupoid), and we see from
  Lemma \ref{L:fib-lem}(iii) that the natural functor $Z\to(Z,H)$ is an
  equivalence of categories over $X$.  It follows that $\nu$ is
  essentially surjective, and is therefore an equivalence of $2$-categories.
\end{proof} 

\begin{para}
  Now suppose that $S$ is a site and $X \to S$ and $Y\to S$ are stacks
  in groupoids.  The fibered category underlying $X$ has a natural
  site structure as follows: a morphism $y \to x$ is a covering if and
  only if its image in $S$ is a covering.  Equivalently, a sieve over
  $x$ is a covering sieve if and only if its image in $S$ is a
  covering sieve.  (To align this with the way one usually thinks of
  sites in a geometric context, one could equally well let the objects
  of the site in question consist of morphisms $\phi: V \to X$, with
  $V$ an object of $S$ thought of as a fibered category over $S$.  It
  is easy to see that this yields a naturally equivalent site.)  Write
  $X^s$ for the resulting site.
\end{para}

\begin{lem}\label{L:stack-trans} 
Let $F: Y \to X$ and $G: X \to S$ be two categories
  fibered in groupoids.  Suppose $S$ is a site and $X$ is given the
  induced site structure.  If $F$ and $G$ are stacks, then $GF$ is a
  stack.
\end{lem}

\begin{proof} Upon replacing $X$ and $Y$ by equivalent fibered
  categories if necessary, we may assume that $F$ and $G$ have
  compatible cleavages over $S$, so that the concept of ``pullback''
  is defined and functorial with respect to $F$.

  We first claim that $Y$ is a prestack over $S$, i.e., that the Hom
  presheaves are sheaves.  Let $z'\to z$ be a covering in $S$ and
  $y_{1},y_{2}$ objects of $Y_{z}$ with images $x_{1}$ and $x_{2}$ in
  $X_{z}$.  Consider the map of presheaves on $S_{/z}$
  $$\iota:\hom_z(y_1,y_2)\to\hom_z(x_1,x_2).$$ Since $X$ is a stack
  over $S$, the target is a sheaf.  We claim that the pullback of
  $\iota$ via any map $\phi:t\to\hom_z(x_1,x_2)$ over a map $t\to z\in
  S$ is a sheaf on $S_{/t}$.  Let $\tilde\phi:x_1|_t\to x_2|_t$ be the
  map corresponding to $\phi$.  Since $Y$ is a category fibered in
  groupoids over $X$, the pullback is a torsor under
  $\hom_{x_1|_t}(y_1|_t,\tilde\phi^{\ast}y_2|_t)$, which is a sheaf
  since $Y$ is a (pre)stack over $X$.  Thus, $\iota$ is a map of
  presheaves with the codomain a sheaf and which is relatively
  representable by sheaves; it follows that the domain is a sheaf,
  which shows that $Y$ is a prestack over $S$.

  Now let $\widebar Y$ be the stackification of $Y$ as a fibered
  $S$-category, constructed in \S II.2 of \cite{giraud} or in Lemme
  3.2 of \cite{l-mb}.  
Since $Y$ is a prestack over $S$, the map of
  fibered categories $Y\to\widetilde Y$ is fully faithful on fiber
  categories (a monomorphism of fibered categories).  The universal
  property of stackification yields a map of $S$-stacks $Q:\widebar Y\to
  X$; letting $\widetilde Y=(\widebar Y,Q)$ (which is still an
  $S$-stack, as it is equivalent to $\widebar Y$), there results a diagram
  $$
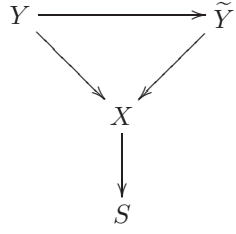
\xymatrix{Y\ar[rr]\ar[rd] & & 
    \widetilde Y\ar[ld]\\
    & X\ar[d] & \\
    & S & \\}$$ of fibered categories over $S$.  Let $z\in S$, $x\in
  X_{z}$, and $t\in\widetilde Y_{x}$.  By the construction of
  $\widebar Y$ (and the equivalence $\widebar Y\cong\widetilde Y$),
  there is a covering $z'\to z$, an object $y'\in Y_{z'}$, and an
  arrow $\sigma:y'\to t|_{z'}$ over the identity of $z'$.  Let
  $\rho:x''\to x|_z'$ be the image of $\sigma$ in $X$, so that $\rho$
  is an arrow over the identity of $z'$, hence has an inverse.  There
  results an isomorphism $(\rho^{-1})^{\ast}y'\to t|_{z'}$ over
  $x|_{Z'}$.  Thus, $Y\to\widetilde Y$ is a morphism of $X$-stacks
  which is fully faithful and such that every object of $\widetilde Y$
  is locally in the image of $Y$.  It follows that it is a
  $1$-isomorphism of stacks, and thus that $Y$ is $1$-isomorphic to a
  stack on $X$ which is a stack over $S$.  The result follows.
\end{proof}

\begin{lem}\label{L:stacky-stack}
  Suppose $X$ and $Y$ are $S$-stacks and $F:Y\to X$ is a functor over
  $S$.  The fibered replacement $(Y,F)$ is a stack on $X$ and on $S$.
\end{lem}
\begin{proof}

  We can describe $(Y,F)$ as the graph of $F$ in the following way:
  the second projection and the functor $F$ yield two maps
  $\pr_2,F\circ\pr_1:Y\times X\to X$.  Given an object $s\in S$ and a
  map $f:s\to Y\times X$, we thus get a sheaf $\isom_s(\pr_2\circ
  f,F\circ\pr_1\circ f)$ on $s$.  It is easy to see that this is a
  ``construction locale'' in the sense of \S{14} of \cite{l-mb}; in
  fact, this defines a sheaf on $(Y\times X)^s$.  It is well-known
  that sheaves $\ms Z$ on $(Y\times X)^s$ are in bijection with
  representable morphisms of $S$-stacks $Z\to Y\times X$.  The $S$-stack
  in question here is easily seen to be $(Y,F)$.  To see that $(Y,F)$
  is an $X$-stack, we note that we have a diagram $(Y,F)\to Y\times
  X\to X$ where each arrow is a stack (the first is even a sheaf).
  Thus, the result follows from Lemma \ref{L:stack-trans}.
\end{proof}

\begin{prop}\label{P:stacky-stack}
  The equivalence of Proposition \ref{P:fib-equiv} restricts to an
  equivalence of $S$-stacks over $X$ and $X^s$-stacks.
\end{prop}
\begin{proof}
  This follows immediately from Lemma \ref{L:stacky-stack}.
\end{proof}

\begin{cor} Suppose $\ms X$ is an algebraic stack.  Given 
  $\alpha \in \H ^ 2 (\ms X _ {fppf}, \m _ n)$, there is a morphism of
  algebraic stacks $\pi: \ms Y \to \ms X$ and a map $\m_n\to\ms I(\ms
  Y)$ such that
  \begin{enumerate}
  \item for any algebraic space $\phi: V \to \ms X$, $\ms Y \times _
    {\ms X} V \to V$ is a $\m _ n$-gerbe representing $\alpha|_{V}$;
  \item the natural morphism of sites associated to $\pi$ is a $\m _
    n$-gerbe representing $\alpha$.
  \end{enumerate}
\end{cor}
\begin{proof}
  This follows immediately from Proposition \ref{P:stacky-stack}, once we note
  that a relative $\m_n$-gerbe over an algebraic stack is itself
  algebraic (as algebraicity is local in the smooth topology).
\end{proof}

This allows us to carry over the formalism of twisted
sheaves (including the Grothendieck existence theorem) to algebraic
stacks without any significant changes.

\begin{rem}  The reader uncomfortable with the techniques sketched 
above can translate the results of this paper into the language of 
Azumaya algebras and their deformations in the \'etale topos of a 
Deligne-Mumford stack, at the expense of knowing a bit more about the 
\'etale cohomology of such stacks.  We will not pursue this point.
\end{rem}

\subsection {Vector bundles on a stack and descent of 
formal extensions}

Let $\ms X$ be a tame Deligne-Mumford stack separated and of finite
type over a locally Noetherian algebraic coarse moduli space $\pi:\ms
X\to X$.  Define a stack $\ms V_{\ms X/X}$ (or simply $\ms V$ when
$\ms X$ is understood) on the big fppf site of $X$ whose fiber
category $\ms V_T$ over $T=\spec A\to X$ is the groupoid of locally
free sheaves on $\ms X\times_X T$.

\begin{prop} With the above notation, the stack $\ms V_{\ms X/X}$ is 
an Artin stack smooth over $X$.
\end{prop}
\begin{proof} In the general case, it is easiest to proceed by Artin's
  representability theorem \cite{artin}.  (This requires descending
  $\ms X$ over an excellent base, which is easily done by the standard
  techniques in \S 8 of \cite{ega4-3}.)  This is made significantly
  easier by the fact that the natural deformation and obstruction
  modules for any object $V\in\ms V_T$ (over \emph{affine} $T=\spec
  A$) with respect to the trivial infinitesimal extension $A+M$ (in
  the notation of \cite{artin}) vanish.  Indeed, these are described
  by the higher cohomology of $\send_{\ms X_{T}}(V)\tensor M$, which
  vanishes because the coarse moduli space $T$ of $\ms X_T$ is affine
  and $\ms X_T$ is tame.  The Schlessinger conditions (condition S1 in
  \cite{artin}) then follow trivially, as the set of deformations of
  an object $a\in\ms V(A)$ over any infinitesimal extension $A'\to A$
  is a singleton (!).  The infinitesimal automorphisms of $V\in\ms
  V(A)$ with respect to the extension $A+M$ are given by
  $\End(V)\tensor_{A} M$, which is a finite $A$-module compatible with
  arbitrary base change.  Finally, the Grothendieck existence theorem
  for $\ms V$ follows immediately from Proposition 2.1 of
  \cite{olsson-starr}.  Putting this all together easily yields the
  result; the reader who is interested in this level of generality is
  encouraged to work out the rest of the details.

  In this paper we will only actually need this result when $\pi$ is a
  flat morphism, and in this case it is easy to give a direct proof:
  it suffices to work \'etale locally on $X$, so we may assume that
  $X=\spec B$ and $\ms X=[\spec C/G]$, with $C$ a finite locally free
  $B$-algebra and $G$ a finite group.  In this case, the stack $\ms V$
  is identified with the substack of $G$-equivariant locally free
  $B$-modules which admit an equivariant $C$-action.  It is trivial
  that locally free $G$-equivariant $B$-modules form an algebraic
  stack locally of finite presentation over $B$, as adding a
  $G$-structure amounts to choosing sections of the automorphism sheaf
  satisfying algebraic conditions coming from a presentation of $G$.
  Since everything is finite and flat, it is also easy to see that the
  addition of a $C$-structure is an algebraic condition (described in
  terms of homomorphisms between tensor products).  The result
  follows.
\end{proof}

\begin{cor}\label{C:formal extension} Let $\ms X$ be a tame separated 
  Deligne-Mumford stack of finite type over a locally Noetherian
  algebraic space with coarse moduli space $\pi:\ms X\to X$.  Suppose
  $X=\spec A$ is affine and $I\subset A$ is an ideal.  Let
  $X'=\spec\widehat A_I$, $U=X\setminus Z(I)$, and $U'=U\times_X X'$.
  Given a locally free sheaf $V$ on the open substack $\ms X\times_X
  U$ and an extension $V'$ of $V_{U'}$ to $\ms X\times_X {X'}$, there
  is a locally free sheaf $W$ on $\ms X$ extending $V$ and such that
  $W_{X'}\cong V'$.
\end{cor}
\begin{proof} This follows from the previous proposition combined with
  Theorem 6.5 of \cite{moret-bailly}.  The history surrounding this
  type of descent problem and various situations in which it holds are
  also described in [{\em ibid.}].
\end{proof}

\section {Splitting ramification}\label{S:splitting on curves}

\subsection{Splitting by cyclic orbifold covers}
\label{sec:splitt-cycl-orbif}

We prove a result in the spirit of Saltman \cite{saltman-fix, saltman}
on moving the ramification in a cyclic cover.  Let $X \to T$ be a
regular surface flat over a local Dedekind scheme with infinite
perfect residue field.  Denote the reduced special fiber by $X _ 0$.
Let $\alpha \in \Br (\kappa(X))'$.  By standard methods \cite{lipman},
we may blow up $X$ and assume that $\ram(\alpha) + X _ 0$ is a strict
normal crossings divisor, say $\ram(\alpha) = D = D _ 1 \cup \dots
\cup D _ r$.

\begin{prop}\label{P:making nth power}  
  Given $n > 0$, there is a $T$-very ample normal crossings divisor
  $C$ such that
  \begin{enumerate}
  \item $C + D$ has normal crossings;
  \item there is an ample divisor $H$ with $C + D \sim n H$.
  \end{enumerate}
\end{prop}

\begin{proof} First, note that $X$ is projective over $T$ (as any
  multisection which intersects every component of the closed fiber is
  $T$-ample). 
  It is easy to see that if $H$ is ``sufficiently very ample'' then
  $C=nH-D$ is very ample, thus proving second statement.

  The crossings of $D$ occur only in the closed fiber.  Since $k$ is
  perfect and infinite, it follows from Bertini's theorem that there
  is a section of $n H - D$ over the reduced special fiber which
  contains only smooth points of $D_0$.  This hyperplane lifts over
  all of $T$, giving rise to $C$.  Since the intersection of $C$ with
  $D$ is reduced, it follows that $C$ has regular components.  Thus,
  $C + D$ has normal crossings.
\end{proof}

Let $X'\to X$ denote the cyclic cover resulting from extracting an
$n$th root of the divisor $C+D$, as in section 2.4 of \cite{kollar-mori}.  By
local calculations, one can see that $X'$ is normal with
finite quotient singularities, so that the local rings at singular
points look like $R[z]/(z^n-xy)$, where $R$ is a local ring of $X$ and
$x$ and $y$ are a system of parameters.

\begin{lem}\label{L:orbifold-res}
  There is a Zariski $A_{n-1}$-orbisurface $\ms X$ with
  coarse moduli space $\pi:\ms X\to X'$ such that $\pi$ is an
  isomorphism over the regular locus of $X'$.
\end{lem}
\begin{proof}
  This is made quite easy by the fact that the singularities of $X'$
  are quotient singularities in the Zariski topology.  Let
  $x_1,\ldots,x_m\in X'$ be the (closed) singular points and let
  $V=X'\setminus\{x_1,\ldots,x_m\}$ be the regular locus.  For each
  $x_i$, there is a regular $T$-scheme $U_i$ with a $\m_n$ action and
  a $\m_n$-invariant morphism $U_i\to X'$ such that $\m_n$ acts freely
  in codimension $1$ and fixes a single closed point.  It follows that
  the induced map from the (regular) quotient stack $\ms
  U_i:=[U_i/\m_n]\to X'$ is an isomorphism over
  $V\times_X'[U_i/\m_n]$.  Gluing the $\ms U_i$ and $V$ together
  yields $\ms X$.
\end{proof}

The following result is the key to the rest of this paper.  By using
an orbifold resolution of the quotient singularities of $X'$ instead
of the standard minimal resolution in the sense of Lipman
\cite{lipman} (as Saltman used), we avoid having to think about
additional ramification divisors.  This results in an optimally efficient
way of splitting the ramification of $\alpha$, in the sense that the
function field extension has minimal degree.

\begin{prop}\label{P:unram-I}
  Let $\ms X$ be the stack constructed in Lemma \ref{L:orbifold-res}.  There
  is an Azumaya algebra $\ms A$ on $\ms X$ with Brauer class
  $\alpha|_{\ms X}$.  In particular, there is a $\m_n$-gerbe $\ms
  G\to\ms X$ representing an extension of $\alpha$ to $\ms X$.
\end{prop}
\begin{proof}
  Let $V$ be the non-stacky locus of $\ms X$ (which is isomorphic via
  the natural map $\ms X\to X'$ to the $V$ in the proof of Lemma 
  \ref{L:orbifold-res}).  It is a standard fact in the ramification
  theory of central simple algebras that $\alpha|_V$ is unramified, as
  $V\to X$ extracts an $n$th root of a generic uniformizer of each
  component of the ramification divisor (see e.g.\ Proposition 1.3 of
  \cite{saltman-fix}).  Thus, there is an Azumaya algebra $\ms
  A^{\circ}$ on $V$ representing $\alpha|_V$.  Taking the reflexive
  hull and using the standard purity theorems on an \'etale cover of
  $\ms X$ shows that $\ms A^{\circ}$ extends to an Azumaya algebra on
  all of $\ms X$, as required.  (In place of \'etale cohomology, one
  could make the following alternative argument: proceeding as above,
  we end up with a sheaf of algebras $\ms A$ on $\ms X$ with the
  property that the natural map of locally free sheaves $\ms
  A\tensor\ms A^{\circ}\to\send(\ms A)$ is an isomorphism in
  codimension $1$.  Since $\ms X$ is regular, this map must be an
  isomorphism, which implies that $\ms A$ is Azumaya.)
\end{proof}

\begin{cor}\label{C:unram-stack-I}
  Given $X\to T$ and $\alpha$ as above, there is an
  $A_{n-1}$-orbisurface, proper and flat over $T$, admitting a map
  $\ms X\to X$ which is finite of generic degree $n$ and such that
  $\alpha|_{K(\ms X)}$ extends to a class in $\Br(\ms X)$.
\end{cor}

\subsection{Splitting by root constructions}
\label{sec:splitt-root-constr}

In this section, we assume that $X$ is a connected regular excellent
scheme of dimension $2$ and $D\subset X$ is a divisor with simple
normal crossings.  Write $U=X\setminus D$.  Suppose $\alpha\in\Br(U)$
is a Brauer class of period $n$ with $n\cdot 1\in\Gamma(X,\ms
O_X)^{\ast}$.  We assume that $\Gamma(X,\ms O_X)$ contains a primitive
$n$th root of unity $\zeta$.

Writing $D=D_1+\cdots+D_r$ as a sum of regular irreducible components,
let $X\{D^{1/n}\}$ denote the stack
$X_{(D_1,\ldots,D_r),(n,\ldots,n)}$ in the notation of Definition
2.2.4 of \cite{cadman}.  We have a natural morphism $X\{D^{1/n}\}\to
X$ which is an isomorphism over $U$ and which is universal among
morphisms for which the pullback of each $D_i$ is a the $n$th power of
an effective Cartier divisor.

\begin{prop}\label{sec:splitt-root-constr-1}
  There is a unique class $\widetilde\alpha\in\Br(X\{D^{1/n}\})$ whose
  restriction to $U$ is $\alpha$.  Moreover, $\widetilde\alpha$ has
  period $n$.
\end{prop}
\begin{proof}
  Uniqueness follows from the fact that the map
  $\Br(X\{D^{1/n}\})\to\Br(k(X))$ is injective.  To prove that
  $\alpha$ extends, we proceed as follows.  First, it suffices to
  prove that it extends in the analogous situation in which $X=\spec
  R$ with $R$ a dvr and $D$ is the closed point of $X$.  Indeed, it
  follows from purity that the functor $V\mapsto\Br(V)$ is a Zariski
  sheaf on $X\{D^{1/n}\}$; if we achieve the result over dvrs, then we
  know that $\alpha$ extends over some open $V\subset X\{D^{1/n}\}$
  containing the generic points of $D$.  Invoking purity again, we
  know that the restriction map $\Br(X\{D^{1/n}\})\to\Br(V)$ is an
  isomorphism, establishing the result.  In less fancy language:
  given an Azumaya algebra $A$ over $V$, its pushforward to
  $X\{D^{1/n}\}$ will be reflexive, hence locally free (as $X$ is a
  regular surface) and Azumaya in codimension $1$, whence it is
  Azumaya since the degeneracy locus of the trace is of codimension
  $1$ (and thus empty in this case).

  So we may suppose $R$ is a dvr with uniformizer $t$ and residue
  field $\kappa$.  The ramification of $\alpha$ is given by some
  cyclic extension $\kappa(\widebar a^{1/n})\kappa$, where $a\in
  R^\times$ is a unit.  Letting $(a,t)_n$ denote the cyclic algebra of
  degree $n$, we see that the class $\alpha-(a,t)_n$ is unramified.
  Thus, to prove the result, it suffices to prove it for the algebra
  $(a,t)_n$.  Note that a local model for $X\{D^{1/n}\}$ is given by
  the stack quotient $[\spec R[t^{1/n}]/\m_n]$.  To show that
  $(a,t)_n$ extends to an unramified class over $X\{D^{1/n}\}$, it
  thus suffices to find a $\m_n$-equivariant Azumaya order in
  $(a,t)_n\tensor K(R)(t^{1/n})$.  

  Recall that $(a,t)_n$ is generated by $x$ and $y$ such that $x^n=a$,
  $y^n=t$, and $xy=\zeta yx$.  We can choose new generators for
  $(a,t)_n\tensor K(R)(t^{1/n})$ by letting $\widetilde y=y/t^{1/n}$.
  The $R[t^{1/n}]$-order generated by $x$ and $\widetilde y$ is then
  Azumaya (as its reduction modulo $t^{1/n}$ is the split central
  simple algebra $(\widebar a,1)_n$), and has a $\m_n$-action arising
  from the action on $t^{1/n}$.
\end{proof}

The reduced preimage of $D$ in $X\{D^{1/n}\}$ is a residual curve of
the type studied in Section \ref{sec:residual} above.  Thus, $\alpha$
gives rise to a Brauer class of period $n$ on this residual curve.
The moduli of uniform twisted sheaves on such a residual curve
developed in Section \ref{sec:moduli-unif-twist} will give us
information that can be deformed off of the ramification locus and
back to the function field of $X$ (in favorable circumstances).

\section{Some cohomology}
\label{S:coho-mama}

In this section, we gather together the cohomological computations
which we will need throughout the rest of the paper.  This includes a
study of the Brauer group and second cohomology of $\m_n$ on the
classifying stack $\B{\m_n}$ as well as a discussion of the second
cohomology of certain stacky surfaces which will play a fundamental
role starting in Section \ref{S:splitting on curves}.

\subsection{Computation of $\Br(\B{\m_{n}})'$}
\label{sec:brauer-group}

\begin{lem}\label{L:no-funnies}
  Let $Y$ be a scheme such that $\Br(Y)=\H^2(Y,\G_m)_{\text{\rm
      tors}}$.  The natural map
  $$\Br(\B{\m_{n,Y}})\to\H^2(\B{\m_{n,Y}},\G_m)_{\text{\rm tors}}$$ is
  an isomorphism.
\end{lem}
\begin{proof}
  The natural map $Y\to\B{\m_{n,Y}}$ is finite locally free of degree
  $n$.  Thus, an element $\alpha\in\H^2(\B{\m_{n,Y}},\G_m)$ lies in
  $\Br(\B{\m_{n,Y}})$ if and only if $\alpha|_Y$ lies in $\Br(Y)$ (see
  e.g.\ 3.1.4.5 of \cite{period-index-paper}).  But this is ensured by the
  hypothesis.
\end{proof}

\begin{remark}
  The lemma notably holds when $Y$ is the spectrum of a field.  By
  Gabber's theorem \cite{dejong-gabber}, we thus see that the Brauer
  group and cohomological Brauer group of $\B{\m_{n,Y}}$ coincide
  whenever $Y$ is a quasi-compact separated scheme admitting an ample
  invertible sheaf.
\end{remark}

\begin{lem}\label{L:bmn-lem}
  If $\widebar L$ is separably closed and $n$ is invertible in
  $\widebar L$ then $\Br(\B{\m_n})'=0$.
\end{lem}
\begin{proof}
  In this case, $\widebar L$ contains a primitive $n$th root of unity,
  so that it is enough to show $\Br(\B{\Z/n\Z})=0$.  Thus, it suffices
  to show that any $\Z/n\Z$-equivariant central simple $\widebar
  L$-algebra $A$ of degree invertible in $\widebar L$ is isomorphic to
  the endomorphism algebra of a $\Z/n\Z$-equivariant vector space $V$.
  Since $\widebar L$ is separably closed, the algebra underlying $A$
  is simply $\M_r(\widebar L)$ for some integer $r$ invertible in $L$.  By
  the Skolem-Noether theorem, the generator of $\Z/n\Z$ acts on
  $\M_r(\widebar L)$ via conjugation by some element
  $\phi\in\GL_r(\widebar L)$ such that $\phi^n$ is a scalar
  $\gamma$.  Since $n\in L^{\times}$, there is some $\eps$ such
  that $\eps^n=\gamma^{-1}$.  Setting $\phi'=\eps\phi$ yields an
  element of order $n$ of $\GL_r(\widebar L)$, i.e., a
  $\Z/n\Z$-equivariant structure on $\widebar L^r$, giving rise to the same
  equivariant structure on $\M_r(\widebar L)$.  The result follows.
\end{proof}

\begin{prop}\label{P:basic-seq-mn}
  Let $L$ be a field.  There is a split exact sequence
$$0\to\Br(L)'\to\Br(\B{\m_n})'\to\H^1(\spec L,\Z/n\Z)\to 0$$
compatible with base extension.
\end{prop}
\begin{proof}
  This comes from the Leray spectral sequence for the projection
  $\pi:\B{\m_n}\to L$, using (1) $\R^1\pi_{\ast}\G_m=\Z/n\Z$ and (2)
  $\R^2\pi_{\ast}\G_m=0$.  The first fact is simply the computation of
  the character group of $\m_n$, while the second follows from Lemma
  \ref{L:bmn-lem} and Lemma \ref{L:no-funnies}.  The splitting comes
  from the natural map $\spec L\to\B{\m_n}$.
\end{proof}

\begin{cor}\label{C:stack-pt-det}
  Let $L$ be a field and suppose $n\in L^{\times}$.  Suppose that there
  is a positive integer $d$ such that 
  every $\alpha\in\Br(L)[n]$ satisfies $\ind(\alpha)|\per(\alpha)^d$.
  Then every $\beta\in\Br(\B{\m_n})[n]$ satisfies
  $\ind(\beta)|\per(\beta)^{d+1}$.
\end{cor}
\begin{proof} Suppose $\beta\in \Br(\B{\m_n})$ has order $m$, where
  $m$ is a positive divisor of $n$.  Using \label{P:basic-seq-mn}, we
  can write $\beta=\beta'+\beta''$ with $\beta'$ in $\Br(L)[m]$ and
  $\beta''\in\H^1(\spec L/\Z/m\Z)\subset\H^1(\spec L,\Z/n\Z)$. By
  assumption, there is a field extension $L'/L$ of degree $m^d$
  annihilating $\beta'$.  As there is clearly a finite \'etale
  extension $L''/L$ of degree $m$ annihilating $\beta''$, we conclude
  from the compatibility of Proposition \ref{P:basic-seq-mn} with
  field extension there is a finite \'etale morphism $e:\spec
  L'''\to\spec L$ of degree $m^{d+1}$ such that
  $\beta|_{\B{\m_n}\tensor_L L'''}=0$.
\end{proof}

\subsection{Computation of $\H^2(\B\m_n,\m_n)$}
\label{sec:cohomology-m_n}

\begin{prop}\label{P:br-m_n}
  Let $L$ be a field.  The Kummer sequence and the natural covering
  $\spec L\to\B{\m_n}$ give rise to a canonical isomorphism
$$\H^2(\B{\m_n},\m_n)\cong\Z/n\Z\oplus\Br(L)[n]\oplus\H^1(\spec
L,\m_n)=\Pic(\B\m_n)\oplus\Br(\B\m_n)$$
in which projection to the first factor corresponds to base change to
$\widebar L$ and inclusion of the first factor corresponds to the
first Chern class map $\Pic(\B\m_n)\to\H^2(\B\m_n,\m_n)$.  This
isomorphism is compatible with base extension.
\end{prop}
\begin{proof}
  This follows readily from Proposition \ref{P:basic-seq-mn}.
\end{proof}

\begin{remark}\label{R:explicit-coboundary}
  Let us make explicit the boundary map
  $\delta:\H^1(\B\m_n,\PGL_n)\to\H^2(\B\m_n,\m_n)$ when the base field
  is separably closed (and $n$ is invertible).  Fix an isomorphism
  $\Z/n\Z\simto\m_n$. The cohomology of $\B\m_n$ is then simply the
  group cohomology of $\Z/n\Z$.  Given a homomorphism
  $\Z/n\Z\to\PGL_n$ (corresponding to a class $[\ms
  A]\in\H^1(\B\m_n,\PGL_n)$), choose a lift $\alpha$ of $[1]$ to
  $\SL_n$ (which is possible as $L$ is separably closed).  Concretely,
  given a $\Z/n\Z$-equivariant structure on $\M_n(L)$, choose a matrix
  $A$ with $A^n=1$ such that conjugation by $A$ defines the action of
  $[1]$ on $\M_n(L)$.  Now, $\det A$ is an $n$th root of $1$, say
  $\zeta$.  Setting $\alpha$ equal to $\zeta^{-1/n}A$ yields a lift of
  $[1]$ to $\SL_n$.  We see that $\alpha^n$ equals scalar
  multiplication by $\zeta$, which is identified with an element of
  $\Z/n\Z$ via the chosen isomorphism.  This element is precisely
  $\delta([\ms A])$.
\end{remark}

\begin{para}
  When $n$ is prime and the base field contains a primitive $n$th root
  of unity, then we can explicitly describe the image of the map
  $$\H^1(\B\m_n,\PGL_n)_{\text{\rm triv}}\to\H^2(\B{\m_n},\m_n),$$ where
  $\H^1(\B\m_n,\PGL_n)_{\text{\rm triv}}$ consists of elements whose
  associated Brauer class is trivial when pulled back along $\spec
  L\to\B\m_n$.  (More generally, such a description is available for
  $\B{\Z/n\Z}$ rather than $\B{\m_n}$.)  In order to avoid confusion,
  we fix a prime $\ell$, invertible in $L$, and assume that $L$
  contains a primitive $\ell$-th root of unity.  Via a choice of
  primitive root, we get an isomorphism $\Z/\ell\Z\simto\m_{\ell}$,
  which we implicitly use in what follows.
\end{para}

\begin{defn}\label{sec:cohomology-m_n-3}
  Let $L$ be a field.  Given a $\G_m$-gerbe $\pi:\ms G\to\B\m_{n,L}$,
  a locally free $\ms G$-twisted sheaf $\ms W$ is \emph{regular\/} if
  for some (and hence any) section $\tau:\B\m_{n}\tensor \widebar
  L\simto\ms G\tensor\widebar L$ of $\pi\tensor \widebar L$, the
  pullback $\tau^{\ast}\ms W$ is the sheaf associated to the regular
  representation of $\m_n$.  A locally free $\ms G$-twisted sheaf $\ms
  W'$ is \emph{totally regular\/} if it has the form $\ms W^{\oplus
    n}$ for some regular locally free $\ms G$-twisted sheaf $\ms W$.
\end{defn}

A twisted sheaf on a $\m_n$-gerbe over $\B\m_n$ will be called regular
if the pushforward to the associated $\G_m$-gerbe is regular.

\begin{lem}\label{L:reg-lem}
  Suppose $\ell$ is prime and invertible in $L$ and $\pi:\ms
  G\to\B\m_\ell$ is a non-trivial $\G_m$-gerbe whose pullback to
  $\spec L$ by the natural map $\spec L\to\B\m_\ell$ is trivial.  If
  $\ms V$ is a $\ms G$-twisted sheaf of rank $\ell$, then $\ms V$ is
  regular.
\end{lem}
\begin{proof}
  Fix an isomorphism $\Z/\ell\Z\simto\m_\ell$.  Using this
  identification and the triviality of $\ms G|_L$, if we given $\ms
  V$, the algebra $\pi_\ast:\send(\ms V)$ is identified with a
  $\Z/\ell\Z$-equivariant structure on $\M_\ell(L)$, which is in turn
  identified with an element of order $\ell$ in $\PGL_\ell$.  This can
  be lifted to a matrix $A\in\GL_\ell(L)$ such that $A^\ell=b$ for
  some scalar $b\in L^\times$ with non-zero image in
  $L^\times/(L^\times)^\ell$. Extending scalars to $L':=L(b^{1/\ell})$, we
  see that we can detect regularity by looking at the eigenspace
  decomposition of $A':=b^{-1/\ell}A$: we have that $(A')^\ell=1$, so
  that $(L')^\ell$ corresponds to a representation of $\m_\ell$.  This
  is simply the pullback of $\ms V$ along a section of $\ms
  G\to\B\m_\ell$; the different choices of the $\ell$th root of $b$
  correspond to the different choices of sections of the gerbe and
  have no effect on regularity, as they only permute the eigenspaces.

  By elementary Galois theory, we know that the action of $A$ on the
  vector space $V=L^\ell$ has the property that $V\tensor L^{\text{\rm
      sep}}$ breaks up as a direct sum of distinct one-dimensional
  eigenspaces for the action of $A\tensor L^{\text{\rm sep}}$.  But
  these are the same as the eigenspaces for
  $A'\tensor_{L'}L^{\text{\rm sep}}$.  This implies that $(L')^\ell$
  is the regular representation of $\m_\ell$, as desired.
\end{proof}

\begin{cor}\label{sec:cohomology-m_n-2}
  Suppose $L$ is finite.  Let $\ms G\to\B\m_{\ell}$ be a
  $\m_{\ell}$-gerbe with cohomology class $(a,0,\beta\neq 0)$.  Any
  two $\ms G$-twisted sheaves of rank $\ell$ are isomorphic.
\end{cor}
\begin{proof}
  Let $\ms V_1$ and $\ms V_2$ be two such sheaves. By Lemma \ref{L:reg-lem},
  both $\ms V_1$ and $\ms V_2$ are regular, so that the scheme
  $\isom(\ms V_1,\ms V_2)$ is a torsor under a torus.  By Lang's
  theorem, it has a rational point.
\end{proof}

\begin{prop}\label{P:detcompm_n}
  Let $\ms A$ be an Azumaya algebra of degree $\ell$ over
  $\B\m_{\ell}$.  Suppose that the Brauer class of $\ms A$ is
  non-trivial but that $\ms A|_{\spec L}=0$. 
Via the isomorphism of Proposition 
  \ref{P:br-m_n}, the class of $\ms A$ in
  $\H^2(\B\m_{\ell},\m_{\ell})$ has the form
  \begin{enumerate}
  \item  $(1,0,\alpha)$ for some
    $\alpha\neq 0$ if $\ell=2$;
  \item $(0,0,\alpha)$ for some
    $\alpha\neq 0$ if $\ell$ is odd.
  \end{enumerate}
  Thus, given a $\m_{\ell}$-gerbe $\ms Y\to\B\m_{\ell}$ with
  cohomology class $(a,0,\beta\neq 0)$, every locally
  free $\ms Y$-twisted sheaf of rank $\ell$ has determinant
  $a-\delta_{2\ell}\in\Z/\ell\Z$, where $\delta$ is the Kr\"onecker
  delta function.
\end{prop}
\begin{proof}
  Let $\ms G\to\B\m_n$ be the $\G_m$-gerbe of trivializations of $\ms
  A$.  We have that $\ms A|_{\ms G}=\send(\ms V)$ for some locally
  free $\ms G$-twisted sheaf of rank $\ell$. By Lemma \ref{L:reg-lem}, we
  know that $\ms V$ is regular.  If we represent $\ms A$ by an element
  $\widebar A$ of order $\ell$ in $\PGL_\ell$, then the regularity
  implies that any lift $A$ of $\widebar A$ in $\GL_{\ell}\tensor
  L^{\text{\rm sep}}$ such that $A^\ell=1$ gives the regular
  representation of $\Z/\ell\Z$ on $(L^{\text{\rm sep}})^\ell$, and
  thus has trivial determinant if $\ell$ is odd, or determinant $[1]$
  if $\ell=2$.
\end{proof}

\subsection{The Brauer group of $\B(\m_\ell\times\m_\ell)$ and Saltman's meteorology}
\label{sec:brau-group-bm_ellt}

In this section, we compute the Brauer group of the classifying stack
$\B(\m_\ell\times\m_\ell)$ and compare the results to Saltman's recent
``meteorology'' of crossings point in the ramification divisor of
Brauer classes on surfaces \cite{saltman-meteorology}.

Let $L$ be a field containing a chosen primitive $\ell$th root of
unity; we will use this to identify $\m_\ell$ with $\Z/\ell\Z$ in what
follows (although we will keep the notation canonical).  Write
$\Br_0(\B(\m_\ell\times\m_\ell))$ for the kernel of the natural map
$\Br(\B(\m_\ell\times\m_\ell))\to\Br(L)$ induced by pullback to the
point $u:\spec L\to\B(\m_\ell\times\m_\ell)$.  The structural morphism
of $\B(\m_\ell\times\m_\ell)$ and the point $u$ induce a natural
splitting
$\Br(\B(\m_\ell\times\m_\ell))\simto\Br(L)\oplus\Br_0(\B(\m_\ell\times\m_\ell))$.

\begin{lem}\label{sec:brau-group-bm_ellt-6}
  The group $\Br_0(\B(\m_\ell\times\m_\ell)$ classifies group scheme extensions
$$1\to\G_m\to G\to\m_\ell\times\m_\ell\to 1.$$
In particular, the functor sending the listed extension to the
$1$-morphism $\B G\to\B(\m_\ell\times\m_\ell)$ together with the
inertial trivialization $\G_m\to\ms I(\B G/\B(\m_\ell\times\m_\ell))$
defines a a bijection between the set of isomorphism classes of
extensions and the set of isomorphism classes of $\G_m$-gerbes with
Brauer classes in $\Br_0(\B(\m_\ell\times\m_\ell))$.
\end{lem}
\begin{proof}
  It is standard that the morphism $\B G\to\B(\m_\ell\times\m_\ell)$
  is a $\G_m$-gerbe.  Pulling back the associated Brauer class to $L$
  is equivalent to pulling back the group extension along the map
  $\B\mathbf 1\to\B(\m_\ell\times\m_\ell)$, where $\mathbf 1$ is the
  singleton group scheme.  But this pullback is just $\G_m$, resulting
  in the trivial $\G_m$-gerbe $\B\G_m$ over $L$.  Thus, the gerbe $\B
  G\to\B(\m_\ell\times\m_\ell)$ has Brauer class in
  $\Br_0(\B(\m_\ell\times\m_\ell)$, as claimed.

  On the other hand, if $\ms G\to\B(\m_\ell\times\m_\ell)$ is a
  $\G_m$-gerbe representing a class in
  $\Br_0(\B(\m_\ell\times\m_\ell))$, then we know that (1) $\ms
  G\to\spec L$ is a gerbe (as it is a composition of gerbes), and (2)
  there is a map $\spec L\to\ms G$ over $L$.  Thus, $\ms G\cong\B{G}$
  for some group scheme $G$ of finite type over $L$.  Moreover, since
  $\ms G$ is a $\G_m$-gerbe over $\B(\m_\ell\times\m_\ell)$, we get a
  natural expression of $G$ as an extension
$$1\to\G_m\to G\to\m_\ell\times\m_\ell\to 1$$ (using the inertial trivialization $\G_m\simto\ms I(\ms G)$).  

To see that the extensions in question are central, note that we can identify sheaves on $\ms G$
with representations of $G$.  In particular, a locally free $\ms
G$-twisted sheaf of (positive) rank $N$ yields a morphism of sequences
$$\xymatrix{1\ar[r] & \G_m\ar[d]_{\id}\ar[r] & G\ar[d]\ar[r] & \m_\ell\times\m_\ell\ar[r]\ar[d] & 1\\
  1\ar[r] & \G_m\ar[r] & \GL_N\ar[r] &\PGL_N\ar[r] & 1.}$$ Tensoring
any locally free $\ms G$-twisted sheaf of positive rank with the
pullback of the sheaf on $\B(\m_\ell\times\m_\ell)$ corresponding to
the regular representation of $\m_\ell\times\m_\ell$, we can find such
a diagram in which the right vertical arrow is injective.  (This
corresponds to tensoring a $\m_\ell\times\m_\ell$-equivariant Azumaya
algebra with the sheaf of endomorphisms of the regular representation,
which is it self a power of the regular representation and thus
faithful.)  It follows that the middle vertical arrow is injective and
then that the extension is central.
\end{proof}

Using the isomorphism
$\m_\ell\times\m_\ell\simto\Z/\ell\Z\times\Z/\ell\Z$ determined by the
chosen primitive $\ell$th root of unity, a projective representation
$\m_\ell\times\m_\ell\to\PGL_N$ is determined by two elements $X$ and
$Y$ in $\GL_N(L)$ such that $X^{\ell}\in L^{\times}$, $Y^{\ell}\in
L^{\times}$ and $[X,Y]=\gamma\in\m_{\ell}(L)$.  Since $X$ and $Y$ are
determined up to elements of $\G_m(L)$, the triple
$(X^\ell,Y^\ell,\gamma)$ yields an element of
$L^{\times}/(L^{\times})^{\ell}\times
L^{\times}/(L^{\times})^{\ell}\times\m_{\ell}$, and it is clear that
this is an invariant of the isomorphim class of the underlying
extension $1\to\G_m\to G\to\m_\ell\times\m_\ell\to 1$, yielding a
set-theoretic morphism
$$i:\Br_0(\B(\m_\ell\times\m_\ell))\to
(L^\times/(L^\times)^\ell)^{\times 2}\times\m_\ell(L).$$  

\begin{prop}\label{sec:brau-group-bm_ellt-1}
  The map $i$ is an isomorphism of abelian groups.
\end{prop}
\begin{proof}
  We first show that $i$ is a map of groups.  Recall that the group
  law on $\Br_0$ is induced by tensor product of
  $\m_\ell\times\m_\ell$-equivariant matrix algebras.  Given two
  actions $\phi_1$, $\phi_2$ of $\m_\ell\times\m_\ell$ on
  $\M_{N_1}(L)$ and $\M_{N_2}(L)$, we compute the sum in $\Br_0$ by
  taking the tensor product representation
  $\M_{N_1}(L)\tensor\M_{N_2}(L)$.  The elements $X$ and $Y$ then act
  diagonally, and we easily see that
  $i(\phi_1\tensor\phi_2)=i(\phi_1)i(\phi_2)$, as desired.

  We now claim that $i$ is surjective.  Let $(A,B,\gamma)$ be a
  triple.  Consider the (possibly commutative!) $L$-algebra
  $(A,B)_{\gamma}$ with generators $x$ and $y$ subject to the
  relations $x^\ell=A$, $y^\ell=B$, and $xy=\gamma yx$.  The elements
  $x$ and $y$ are units, and their action on $(A,B)_\gamma$ by
  conjugation gives a projective representation $\phi$ of
  $\m_\ell\times\m_\ell$ (by looking at the induced action on the
  endomorphism algebra of the vector space underlying $(A,B)_\gamma$)
  such that $i(\phi)=(A,B,\gamma)$.

  Finally, we claim that $i$ is injective. Suppose given an action of
  $\m_\ell\times\m_\ell$ on $\M_N(L)$ whose invariants are $(1,1,1)$,
  we get commuting elements $X,Y\in\GL_N(L)$ such that $X^\ell=\id$
  and $Y^\ell=\id$.  This gives a representation
  $\m_\ell\times\m_\ell\to\GL_N$ whose image in $\PGL_N$ is the given
  action, which shows that the equivariant matrix algebra is the
  endomorphisms of an equivariant vector space, thus trivializing the
  Brauer class.
\end{proof}
\begin{cor}
  Every class in $\Br_0(\B(\m_\ell\times\m_\ell))$ has period $\ell$.
\end{cor}

\begin{para}
  Of utmost importance will be various indices of the classes in
  $\Br_0$.  The natural point $\spec L\to\B\m_\ell$ yields two maps
  $r_i:\B\m_\ell\to\B(\m_\ell\times\m_\ell)$, $i=1,2$, which we will
  call the restriction maps.  Throughout this section, we fix a
  $\G_m$-gerbe $\ms G\to\B(\m_\ell\times\m_\ell)$ in
  $\Br_0(\B(\m_\ell\times\m_\ell))$ with invariant $i(\ms G)=(A,B,\gamma)$.  We
  will abusively write $r_i:\ms G_{\B\m_\ell}\to\ms G$ for the map
induced by $r_i$; these maps pull back twisted sheaves to twisted
sheaves.
\end{para}
\begin{defn}\label{sec:brau-group-bm_ellt-2}
  A $\ms G$-twisted sheaf $V$ is \emph{totally biregular\/} if
  $r_i^{\ast}V$ is totally regular for $i=1,2$.
\end{defn}

\begin{defn}\label{sec:brau-group-bm_ellt-3}
  The minimal rank of a non-zero totally biregular $\ms G$-twisted
  sheaf will be called the \emph{global index\/}.
\end{defn}

\begin{prop}\label{sec:brau-group-bm_ellt-4}
  Suppose the field $L$ is finite.  
  \begin{enumerate}
  \item If $\gamma\neq 1$ then both the index and global index of
    $\ms G$ are $\ell$.
  \item If $\gamma=1$ but $A\neq 1$ and $B\neq 1$ then both the
    index and global index of $\ms G$ are $\ell$.
  \item If $\gamma=1$ and $A=1$ but $B\neq 1$, then the index of $\ms
    G$ is $\ell$, while the global index of $\ms G$ is $\ell^2$.
  \end{enumerate}
\end{prop}
\begin{proof}
  For the first statement, note that since $\gamma\neq 0$, the cyclic
  algebra $(A,B)_\gamma$ is a central simple algebra of index $\ell$.
  Since $L$ is finite, there is an isomorphism
  $(A,B)_\gamma\simto\M_{\ell}(L)$, thus yielding a projective
  representation $\phi:\m_\ell\times\m_\ell\to\PGL_{\ell}$ with
  $i(\phi)=(A,B,\gamma)$.  This shows that the index is $\ell$; in
  fact, it realizes the given $\m_\ell\times\m_\ell$-action on
  $M_\ell(L)$ as the endomorphisms of some $\ms G$-twisted sheaf $V$
  of rank $\ell$.  Moreover, the restriction of the Brauer class via
  the $r_i$ is given by restricting $\phi$ to each factor $\m_\ell$.
  To see that the sheaf $V$ is totally biregular, we can use the
  explicit description of $(A,B)_{\gamma}$ as follows.  Extending
  scalars to $\widebar L$ and replacing $x$ by $x/A^{1/\ell}$ and $y$
  by $y/B^{1/\ell}$ yields an isomorphism $(A,B)_\gamma\tensor\widebar
  L\simto(1,1)_\gamma\tensor\widebar L$ of
  $\m_\ell\times\m_\ell$-equivariant $\widebar L$-algebras which have
  the same projective representation.  On the other hand, we can realize $(1,1)_\gamma$ explicitly using the elements 
$$X=\begin{pmatrix}
1 & 0 & \cdots & 0\\
0 & \gamma & \cdots & 0\\
\vdots & \vdots & \ddots & \vdots\\
0 & 0 & \cdots & \gamma^{\ell-1}\\
\end{pmatrix}
$$
and
$$Y=\begin{pmatrix}
0&0&\cdots&0&1\\
1&0&\cdots&0&0\\
0&1&\cdots&0&0\\
\vdots&\vdots&\ddots&\ddots&\vdots\\
0&0&\cdots&1&0\\
\end{pmatrix}
$$
from which it is obvious that each of $X$ and $Y$ acts on $L^{\ell}$
as the regular representation of $\m_\ell$.  This shows that the
twisted sheaf $V$ is totally biregular, as desired.

  To prove the second statement, note that since both $A$ and $B$ are
  non-zero, they generate $L^\times/(L^\times)^\ell\cong\Z/\ell\Z$, so
  that there is some $n$ relatively prime to $\ell$ with $A^m=B\in
  L^{\times}/(L^\times)^\ell$.  Consider the field extension
  $L':=L[x]/(x^\ell-A)$.  Scalar multiplication by $x$ and by $x^n$
  yield an $L$-linear action on $L'$ corresponding to a map
  $G\to\GL_{\ell}$; the image in $\PGL_{\ell}$ will be a projective
  representation with invariant $(A,B,1)$.  To check biregularity we
  can again work over $\widebar L$ and replace $x$ by $X/A^{1/\ell}$,
  whence we see that the resulting map $\m_\ell\to\GL_\ell$ is now the
  regular representation.  Thus, the global index is $\ell$, and this
  implies that the index is $\ell$.

  To see the last statement, note that the algebra
  $L[x,y]/(x^\ell-1,y^\ell-B)$ admits an action as in the previous
  paragraph, yielding a $\ms G$-twisted sheaf of rank $\ell^2$ that is
  totally biregular.  On the other hand, since the action of the first
  factor lifts to an element of $\GL_{\ell^2}$, the restriction of
  $\ms G$ via $r_1$ is trivial; since
  $\B\m_\ell\to\B\m_\ell\times\B\m_\ell$ has degree $\ell$, this shows
  that $\ms G$ has index $\ell$.  If the global index were $\ell$,
  there would be a pair of commuting elements $X$ and $Y$ in
  $GL_\ell(L)$ such that $X^\ell=1$ and $Y^\ell=B$.
  By the double centralizer theorem (applied to the subalgebra
  generated by $Y$), $X$ would have to lie in the subalgebra of
  $\M_\ell(L)$ generated by $Y$ (i.e., this algebra is its own
  centralizer). But the only $\ell$th roots of unity in this
  subalgebra are those in $L$, which means that the action of $X$ is
  not (geometrically) conjugate to the regular representation.
\end{proof}

\begin{cor}\label{sec:brau-group-bm_ellt-5}
  Suppose the triple of $\ms G$ is $(A,B,1)$.  Then there is a locally
  free $\ms G$-twisted sheaf $V$ of rank $\ell^2$ such that for any
  invertible $\ms G\tensor\widebar L$-twisted sheaf $\Lambda$, the
  sheaf $V_{\widebar L}\tensor\Lambda^{\vee}$ is the sheaf on
  $\B(\m_\ell\times\m_\ell)$ associated to the regular representation
  of $\m_\ell\times\m_\ell$.
\end{cor}
\begin{proof}
  Since $\gamma=1$, we see that the cohomology class of $\ms G$ has
  the form $\pr_1^\ast\alpha+\pr_2^\ast\beta$, where $\alpha$ is the
  class in $\Br(\B\m_\ell)$ associated to $A$ and $\beta$ is the class
  associated to $B$.  Letting $W_1$ be a regular $\alpha$-twisted
  sheaf and $W_2$ a regular $\beta$-twisted sheaf, we see that
  $\pr_1^\ast W_1\tensor\pr_2^\ast W_2$ is a $\ms G$-twisted sheaf
  with the requisite properties.
\end{proof}

\begin{para}
  The reader familiar with \cite{saltman-meteorology} might notice a
  resemblance between Proposition \ref{sec:brau-group-bm_ellt-4} and
  Saltman's ``meteorology'' of points on ramification divisors for
  Brauer classes on surfaces.  In this paragraph we will explain the
  connection (at least in the special case in which we have proved it;
  it is easy to extend the method above to the general case and
  achieve total parity with Saltman's setup).  We will be content to
  simply sketch the correspondence, as it is a digression from our
  primary purpose here.
\end{para}
Let us briefly remind the reader of Saltman's terminology.  Suppose
$R$ is a regular local ring of dimension $2$ with residue
field $\kappa$, and suppose $s$ and $t$ are generators of the maximal
ideal of $R$.  For any domain $B$ arising in the following, we will
write $K(B)$ for its fraction field.  

Given a Brauer class $\alpha$ over $K$ of period $\ell$ invertible in
$R$ which is ramified only over the (transversely intersecting)
divisors cut out by $s$ and $t$, the ramification theory of Brauer
classes (reviewed in Section \ref{S:splitting on curves} below)
produces two ramification extensions $L_1/K(R/sR)$ and $L_2/K(R/tR)$
of degree $\ell$ which have the same order of ramification over the closed point
$\spec\kappa$ lying in both $\spec R/sR$ and $\spec R/tR$.

There are three possibilities for the configuration of $L_1$ and $L_2$
(up to the symmetry of $s$ and $t$ in the notation):
\begin{enumerate}
\item both $L_1$ and $L_2$ ramify over $\spec\kappa$;
\item both $L_1$ and $L_2$ are unramified over $\spec\kappa$ and are non-trivial;
\item $L_1$ is trivial and $L_2$ is non-trivial and unramified over
  $\spec\kappa$.
\end{enumerate}
In Saltman's terminology, the first case is called ``cold,'' the second case is called ``chilly'', and the third is called ``hot''.

To connect this with the Brauer group of $\B(\m_\ell\times\m_\ell)$,
we give a preview of results which will be fully explained in Section
\ref{S:splitting on curves}.  (The reader may find it easier to read
this after reading that section.)  Let $\ms X\to\spec R$ be the stack
arising by taking the $\ell$th roots of the divisors $\spec R/sR$ and
$\spec R/tR$ in $\spec R$ (see Section \ref{sec:splitt-root-constr}
for definitions).  By \ref{sec:splitt-root-constr-1}, the class
$\alpha$ extends to a class in $\Br(\ms X)$.  In fact, we have an
isomorphism between $\Br(\spec R\setminus(\spec R/sR\cup\spec
R/tR))[\ell]$ and $\Br(\ms X)[\ell]$.  Now $\ms X$ contains
$\B(\m_\ell\times\m_\ell)$ as a closed substack (the reduced stack
structure on the preimage of $\spec\kappa$ under the natural map $\ms
X\to\spec R$).  Reduction defines a morphism
$\Br(\ms X)[\ell]\simto\Br(\B(\m_\ell\times\m_\ell))[\ell]$; when the
residue field $\kappa$ is finite, this finally results in a morphism
$$\Br(\spec R\setminus(\spec R/sR\cup\spec R/tR))[\ell]\simto\Br(\ms X)[\ell]\to\Br(\B(\m_\ell\times\m_\ell))=\Br_0(\B(\m_\ell\times\m_\ell)).$$
If $R$ is complete, the composition is in fact an isomorphism and
$\Br(\ms X)$ is itself $\ell$-torsion (so the left side becomes
$\Br(\ms X)$).

Via this map, the three ramification situations described above (cold,
chilly, and hot) correspond precisely to the three types of Brauer
classes described in Proposition \ref{sec:brau-group-bm_ellt-4}.  It
is interesting to note that ramification information for classes over
the complicated field $K$ gets translated into information about
Brauer classes over the classifying stack $\B(\m_\ell\times\m_\ell)$
that only sees the finite residue field.  The techniques explained
here are a demonstration that the geometry
of stacks makes this connection less silly than it might at first
seem, as one can import flexible global geometric methods into
ramified situations.

\section{Residual curves and moduli of uniform twisted sheaves}
\label{sec:residual}

Let $k$ be a field.  Fix an odd prime $\ell$.

\begin{defn}\label{sec:resid-curv-moduli-1}
  A stack $\ms C\to k$ is a \emph{residual curve\/} if
  \begin{enumerate}
  \item $\ms C$ is a connected tame Deligne-Mumford stack of dimension
    $1$ whose coarse moduli space space $C$ is a simple normal
    crossing curve;
  \item if $C_1,\ldots,C_r$ are the (smooth) irreducible components of
    $C$, then there are Zariski-locally trivial $\m_\ell$-gerbes $\ms
    C'_i\to C_i$ such that $\ms C=\ms C'_1\times_C\ms
    C'_2\times_C\cdots\times_C\ms C'_r$.
  \end{enumerate}
\end{defn}

The curves $\ms C_i:=\ms C\times_C C_i$ are then smooth
Deligne-Mumford stacks of dimension $1$ with generic stabilizer
$\m_\ell$ and a finite closed substack over which the stabilizer is
$\m_\ell\times\m_\ell$.  We will write $\ms C_i^\circ$ for $\ms
C_i\setminus\cup_{j\neq i}\ms C_j$; it is the intersection of the
smooth locus of $\ms C$ with $\ms C_i$ and is a $\m_\ell$-gerbe over a
the dense open subscheme $C_i\cap C^{\text{\rm sm}}\subset C_i$.

The way such curves usually arise is as the stacky loci of the root
construction \cite{cadman} applied to a simple normal crossing divisor
in a smooth surface, or to a simple normal crossing divisor in a fiber
of an arithmetic surface.  That is certainly how they will arise for
us.  In this section we fix a residual curve $\ms C$, with irreducible
components $\ms C_1,\ldots,\ms C_r$.  Write $\kappa_i$ for $\Gamma(\ms
C_i,\ms O_{C_i})$; we thus have that $\ms C_i$ is a geometrically
integral smooth curve over $\kappa_i$.

Fix a $\m_\ell$-gerbe $\ms G\to\ms C$.  Write $\ms G_i:=\ms
G\times_{\ms C}\ms C_i$.  We will write $\pi$ for the natural map $\ms
G\to\ms C$ and, by abuse of notation, for the natural map $\ms
G_i\to\ms C_i$.  Finally, we will write $\ms G_i^\circ$ for $\ms
G\times_{\ms C}\ms C_i^\circ$.

We assume that for each geometric point $c\to C^{\text{\rm sm}}$, the
restriction $\ms G\tensor\kappa(c)$ is trivial in $\H^2(\ms
C\tensor\kappa(c),\m_\ell)$.  The following lemma shows that this is a
harmless assumption if one is concerned only with the Brauer class of
$[\ms G]$.

\begin{lem}\label{sec:resid-curv-moduli-6}
  Given a class $\alpha\in\Br(\ms C)[\ell]$, there is a
  $\m_\ell$-gerbe $\ms G\to\ms C$ such that for every geometric point
  $c\to C^{\text{\rm }}$ the restriction $\ms G|_{\ms
    C\tensor\kappa(c)}$ has trivial class in $\H^2(\ms
  C\tensor\kappa(c),\m_\ell)$.
\end{lem}
\begin{proof}
  Given any collection of characters $\chi_i^{t_i}$ of the generic
  geometric stabilizers $\B\m_\ell$ on $\ms C_i^\circ$, there is an
  invertible sheaf $L$ on $\ms C$ such that for any geometric residual
  gerbe $x:\B\m_{\ell,\widebar\kappa}\to\ms C_i^\circ$ the character
  associated to $x^\ast L$ is $\chi_i^{t_i}$.  Indeed, the existence
  of an invertible sheaf $L_i'$ on $\ms C_i^\circ$ with character
  $\chi_i^{j_i}$ follows from the assumption that $\ms C_i^\circ\to
  C_i^\circ$ is a Zariski-locally trivial $\m_\ell$-gerbe.  A local
  calculation then shows that the invertible sheaf
  $L_i:=L_i'(-\sum_{j\neq i} t_j\ms C_j\cap\ms C_i)$ has the property
  that at a residual gerbe $\xi=\B\m_\ell\times\B\m_\ell\inj\ms
  C_i\cap\ms C_j$, we have $L_i|_{\xi}$ as character
  $\chi_i^{t_i}\tensor\chi_j^{t_j}$.  Gluing the $L_i$ together yields
  $L$, as desired.

  Now let $c_i\to C_i^\circ$ be a set of geometric points, and let the
  image of $\ms G\tensor\kappa(c_i)$ be $n_i\in\H^2(\ms
  C_i\tensor\kappa(c),\m_\ell)$.  Let $L$ be the invertible sheaf as
  in the previous paragraph with characters $\chi_i^{n_i}$.  Replacing
  $\ms G$ by a $\m_\ell$-gerbe representing the class $[\ms
  G]-[L]^{1/\ell}$ (where $[L]^{1/\ell}$ denotes the image of $L$ via
  the Kummer map $\H^1(\ms C,\G_m)\to\H^2(\ms C,\m_\ell)$) yields a
  gerbe $\ms G'$ with the same Brauer class and whose cohomology class
  vanishes at all geometric points in $C^{\text{\rm sm}}$, as desired..  
\end{proof}

In this section, we will
describe certain stacks of $\ms G$-twisted sheaves.

\begin{defn}\label{sec:resid-curv-moduli-2}
  Given an element $a\in\Z/\ell\Z$, a sheaf $\ms F$ on $\ms C_i^\circ$
  will be called \emph{isotypic of type $a$\/} if the geometric
  generic fiber $\ms F\tensor\widebar{\kappa(C_i)}$ has the form
  $(\chi^a)^{\oplus N}$ for some $N$ as sheaves on
  $\B\m_{\ell,\widebar{\kappa(C_i)}}$, where $\chi:\m_\ell\to\G_m$ is
  the natural inclusion character.
\end{defn}

\begin{defn}\label{sec:resid-curv-moduli}
  Given a locally free $\ms G_i$-twisted sheaf $V$, an
  \emph{eigendecomposition\/} of $V$ is a direct decomposition
  $V=V_0\oplus V_1\oplus\cdots\oplus V_{\ell-1}$ such that for all $s$
  and $t$, the sheaf $\pi_\ast\shom(V_s,V_t)|_{\ms C_i^\circ}$ is
  isotypic of type $[t-s]$.
\end{defn}

\begin{lem}\label{sec:resid-curv-moduli-3}
  Let $V$ be a locally free $\ms G_i$-twisted sheaf.  Choose an
  algebraic closure $\widebar k$ of $k$.
  \begin{enumerate}
  \item There is at most one eigendecomposition of $V\tensor\widebar
    k$, up to isomorphism.
  \item If $V\tensor\widebar k$ admits an eigendecomposition, then the
    set of numbers $\{\deg\shom(V_s,V_t)\}$ is independent of the
    choice of $\widebar k$ and of eigendecomposition over $\widebar
    k$.
  \end{enumerate}
\end{lem}
\begin{proof}
  We may assume $k=\widebar k$.  Suppose $V=V_0\oplus\cdots\oplus
  V_{\ell-1}=W_0\oplus\cdots\oplus W_{\ell-1}$ are two
  eigendecompositions of $V$.  Write $\widebar K$ for an algebraic
  closure of the function field of $C_i$.  Let $M$ be an invertible
  $\ms G\tensor\widebar K$-twisted sheaf.  Restricting the two
  decompositions to the geometric generic point yields two
  eigendecompositions $\bigoplus
  (V_i\tensor M^\vee)=\bigoplus (W_i\tensor M^{\vee})$ of the
  $\widebar K$-linear $\m_\ell$-representation $V\tensor\widebar K$.  It follows
  from the standard representation theory of $\m_\ell$ that there is a
 (cyclic) permutation $i\mapsto j(i)$ such that $V_i\tensor\widebar
 K=W_{j(i)}\tensor\widebar K$.  Re-indexing, we may assume that
 $j(i)=i$.  We conclude by faithful flatness that $V_i\tensor
 K=W_i\tensor K$ as subsheaves of $V\tensor K$.

 Now consider the quotient $V\to V_i$.  For $j\neq i$, the fact that
 $V_j\tensor K=W_j\tensor K$ implies that the natural surjection
 $\bigoplus W_i\surj V_i$ factors through a surjection $W_i\surj
 V_i$.  Since $W_i$ and $V_i$ are of finite rank, it follows that this
 surjection is an isomorphism.
\end{proof}

The set of degrees $\{\deg\shom(V_s,V_t)\}$ will be denoted $d(V)$.
The set $d(V)$ always contains $0$ and is symmetric with respect to
negation.  Similarly, the set $r(V):=\{\rk\shom(V_s,V_t)\}$ is
independent of the geometric eigendecomposition.

\begin{defn}\label{sec:resid-curv-moduli-4}
  A locally free $\ms G$-twisted sheaf $V$ is \emph{uniform of
    determinant $\ms O$\/} if 
  \begin{enumerate}
  \item $\det V\cong\ms O$;
  \item for all $i$, if $V|_{\ms G_i}$ admits a geometric
    eigendecomposition then $d(V)=\{0\}$ and $r(V)$ is a singleton.
  \end{enumerate}
\end{defn}

In particular, we have the following lemma.

\begin{lem}\label{sec:resid-curv-moduli-5}
  If $V$ is a uniform $\ms G$-twisted sheaf of rank $\ell^2$ then 
  \begin{enumerate}
  \item for every residual gerbe $\xi:\B\m_{\ell,\kappa}\inj\ms
    C^\circ$ in the smooth locus of $\ms C$, the restriction
    $V|_{\xi}$ is totally regular;
  \item for each residual gerbe $\xi'$ of $\ms C$ of the form
    $\B\m_\ell\times\B\m_\ell$, the restriction $V|_{\xi'}$ is totally biregular.
  \end{enumerate}
\end{lem}
\begin{proof}
  It suffices to prove this after passing to an algebraically closed
  base field, so we may assume that $\kappa$ is algebraically closed
  and therefore that $V|_{\xi}$ has an eigendecomposition
  $V|_{\xi}=V_0\oplus\cdots\oplus V_{\ell-1}$.  Choosing an invertible
  $\ms G|_{\xi}$-twisted sheaf $L$, the sheaf $V|_{\xi}\tensor L^\vee$
  is thus identified with a representation $W$ of $\m_\ell$ admitting
  a direct sum decomposition $W=W_0\oplus\cdots\oplus W_{\ell-1}$ such
  that $\shom(W_s,W_t)$ is of constant rank and is isotypic of weight
  $t-s$ for all $s$ and $t$.  It follows that each $W_i$ is isotypic
  of a fixed rank $r$, from which it follows that
  $V|_{\xi}\cong\rho^{\oplus r}$, giving the first statement.  The
  second statement follows from the first: the two factors of the
  inertia stack $\m_\ell\times\m_\ell$ of $\xi'$ are limits of the
  inertia stacks $\m_\ell$ on the two branches of $\ms C$ meeting at
  $\xi'$.  Since restriction to a factor corresponds to deformation
  into the smooth locus of the corresponding branch, we see that total
  biregularity implies total regularity on the branches.
\end{proof}

\subsection{Moduli of uniform twisted sheaves}
\label{sec:moduli-unif-twist}

Write $\xi_1,\ldots,\xi_m$ for the set of singular residual gerbes of
$\ms C$; thus, each $\xi_i$ is isomorphic to
$\B\m_\ell\times\B\m_\ell$.  Write $\ms G(i)$ for the restriction of
$\ms G$ to $\xi_i$.  By Corollary \ref{sec:brau-group-bm_ellt-5}, for each
$i$ there is a totally biregular $\ms G(i)$-twisted sheaf $V(i)$ which
is geometrically isomorphic to a twist of the regular representation
of $\m_\ell\times\m_\ell$.  We
fix such a choice for each $i=1,\ldots,m$ in this section.

\begin{defn}
  Let $\ms U$ denote the $k$-stack whose objects over $T\to\spec k$
  are locally free $\ms G_T$-twisted sheaves $\ms V$ of rank $\ell^2$
  together with an isomorphism $\det\ms V\simto\ms O$ such that for
  each geometric point $t\to T$,
  \begin{enumerate}
  \item the fiber $\ms V_t$ is uniform as a $\ms G_t$-twisted sheaf;
  \item for each $i$, the restriction of $\ms V$ to $(\xi_i)_t$ is
    isomorphic to $V(i)_t$.
  \end{enumerate}
\end{defn}

The primary goal of this section is to prove the following theorem.

\begin{thm}\label{sec:moduli-unif-twist-1}
  The stack $\ms U$ is a smooth Artin stack over $k$ with quasi-affine
  stabilizers.  If $\ell$ is odd, $\ms U$ is geometrically connected.
\end{thm}

The primary use of the theorem in this paper will be the following corollary.

\begin{cor}\label{sec:moduli-unif-twist-2}
  Suppose $\ell$ is odd.  If $\ms C$ is a residual curve over a finite
  field $L$ and $\ms G\to\ms C$ is a $\m_\ell$-gerbe, then for some
  $N>1$ relatively prime to $\ell$, there is a locally free $\ms
  G$-twisted sheaf of rank $N\ell^2$ and trivial determinant.
\end{cor}
\begin{proof}
  Since $\ms U$ is integral with quasi-affine stabilizers, there is a
  dense open substack $\ms V\subset\ms U$ of the form $[P/\GL_n]$ with
  $P$ an algebraic space (Proposition 3.5.9 of \cite{kresch}).  Since
  $\ms U$ is geometrically integral, so is $\ms V$ and thus $P$.  By
  the Lang-Weil estimates, $P$ has a point over an extension $L/k$ of
  degree prime to $\ell$.  Taking the image of this point in $\ms U$
  yields a locally free $\ms G\tensor_k L$-twisted sheaf of rank
  $\ell^2$ and trivial determinant.  The proof follows from the
  following Lemma.
\end{proof}

\begin{lem}\label{sec:moduli-unif-twist-3}
  Let $X$ be a proper stack over a field $k$ and $L/k$ a finite \'etale algebra
  of degree $d$.  Write $\pi:X\tensor L\to X$ for the
  natural projection.  If $V$ is a locally free sheaf of rank $r$ on
  $X\tensor L$ of determinant $\ms O$ then $\pi_\ast V$ is locally
  free on $X$ of rank $r$ and determinant $\ms O$.
\end{lem}
\begin{proof}
  To check that the determinant is trivial, it suffices to work over
  $\widebar k$ (as $X$ is proper), and thus we may assume that $\spec
  L=\sqcup_{i=1}^e\spec k=\cup p_i$ is a finite product of copies of
  $k$.  Now $\pi_\ast V=\oplus V|_{p_i}$, so that $\det\pi_\ast
  V=\bigotimes\det V|_{p_i}$.  If $\det V\cong\ms O$, then for each
  $i$ we have $\det V|_{p_i}\cong\ms O$, giving the result.
\end{proof}

\begin{para}
  Before proving Theorem \ref{sec:moduli-unif-twist-1}, we first describe a
  geometric avatar of the ramification extensions.
Consider the stack $\ms P_i:=\ms Pic^{(1)}_{\ms G^\circ_i/C_i^\circ}$
parametrizing $C_i^\circ$-flat families of $\ms G_i^\circ$-twisted
sheaves on $\ms C_i^\circ$.  Standard methods show that $\ms P_i\to
C_i^\circ$ is a $\G_m$-gerbe over a  $R_i\to C_i^\circ$.
\end{para}
\begin{prop}\label{sec:moduli-unif-twist-4}
  With the preceding notation, there is a cyclic \'etale covering
  $R_i\to C_i^\circ$ such that $\ms P_i$ is a $\G_m$-gerbe over $R_i$.
\end{prop}
Intrinsically, $R_i$ is the rigidification of $\ms P_i$ with respect
to $\G_m$, in the notation of \cite{rigid}.  
\begin{proof}
  It is standard that $\ms P_i$ is a $\G_m$-gerbe over an algebraic
  space $R_i$.  Moreover, elementary deformation theory shows
  that $\ms P_i\to C_i^\circ$ is smooth, so that $R_i\to C_i^\circ$ is
  also smooth.  
On the other hand, the relative Picard stack of
  $\ms C_i^\circ$ over $C_i^\circ$ is a $\G_m$-gerbe over the constant
  group scheme $\Z/\ell\Z$.  Tensoring defines an action
  $\Z/\ell\Z\times R_i\to R_i$ over $C_i^\circ$ which is simply
  transitive on geometric fibers.  We conclude that $R_i\to C_i^\circ$
  is naturally a $\Z/\ell\Z$-torsor, as desired.
\end{proof}

\begin{remark}
  It is not too difficult to see that $R_i$ is precisely the
  ramification extension of $\alpha$ at $C_i$.  The easiest proof goes
  by noting that they split one another.
\end{remark}

There is an elementary criterion for the existence of an
eigendecomposition.

\begin{prop}\label{sec:moduli-unif-twist-5}
  Given a field extension $L/\kappa_i$, a non-zero locally free $\ms
  G_i$-twisted sheaf $V$ has an eigendecomposition if and only if
  $R_i\tensor_{\kappa_i}L$ is disconnected.
\end{prop}
\begin{proof}
  It suffices to assume $\kappa_i=L$ and to work at the generic point
  $\spec K$ of $C_i$.  We first claim that there exists a $V$ with an
  eigendecomposition if and only if the Brauer class of $\ms G$ is
  pulled back from $\spec K$.

  To see this, suppose $V=V_0\oplus\cdots\oplus V_{\ell-1}$ is an
  eigendecomposition.  Then $\shom(V_0,V_0)$ is an Azumaya algebra on
  $\B\m_{\ell,K}$ with the same Brauer class as $\ms G$ and on which
  $\m_\ell$ acts trivially.  It follows that the Brauer class of $\ms
  G$ is pulled back from $\spec K$.  On the other hand, if the class
  is a pullback, then let $A$ be an Azumaya algebra over $\spec K$
  representing $[\ms G]$.  The category of $\ms G$-twisted sheaves is
  naturally equivalent to the category of $\m_\ell$-equivariant
  right $A$-modules.  But the representation theory of $\m_\ell$ on
  vector spaces over a division ring is the same as that of a field,
  so we see that any module admits an eigendecomposition.

  Thus, to prove the Proposition it suffices to prove that $R_i$ is
  split if and only if $[\ms G]$ is pulled back from $\spec
  K$. Suppose $R_i$ is split, and let $s:\spec K\to R_i$ be any
  section.  By definition, the obstruction to lifting $s$ to an
  invertible $\ms G$-twisted sheaf vanishes if and only if $[\ms G]$
  vanishes.  This obstruction is thus (at the very least) an element
  of $\Br(K)$ whose pullback under the injective map
  $\Br(K)\inj\Br(\B\m_{\ell,K})$ generates the same cyclic subgroup as
  $[\ms G]$.  It follows that $[\ms G]$ is a pullback.  Conversely,
  assume that $\ms G$ is a pullback from $K$, say $\ms G=\widebar{\ms
    G}\times_{\spec K}\B\m_{\ell,K}$.  The $\G_m$-gerbe
  associated to $\widebar{\ms G}$ over $\spec K$ is naturally isomorphic to the
  relative twisted Picard stack of ${\ms G}$ over $K$, so that
  $\Pic_{\widebar{\ms G}/\spec K}=\spec K$.
  Given a $K$-scheme $T$, it is easy to see any $\ms G_T$-twisted
  sheaf $L$ has the form $L\tensor\Lambda$ with $L$ an invertible
  $\widebar{\ms G}_T$-twisted sheaf and $\Lambda$ an invertible sheaf
  on $\B\m_{\ell,T}$.  This shows that the natural map
  $\Pic_{\B\m_\ell/\spec K}\times\Pic_{\widebar{\ms G}/\spec
    K}\to\Pic_{\ms G/\spec K}$ is an isomorphism, which shows that
  $R_i$ is split, as desired.
\end{proof}

\begin{para}
  We now prove Theorem \ref{sec:moduli-unif-twist-1}.  Since the
  statement is geometric, it suffices to work over the algebraic
  closure of the base field.  Thus, \emph{we will assume that the base
    field $k$ is algebraically closed for the rest of this section\/}.
  Without loss of generality, we may assume that $i=1,\ldots, s$ are
  the indices such that the coverings $R_i\to C_i^\circ$ are split.
  For each $i=1,\ldots,s$, we fix an invertible $\ms G_i$-twisted
  sheaf $\Lambda_i$.  
\end{para}
The idea of the proof is the following. Basic deformation theory shows
that $\ms U$ is smooth, so it suffices to show that it is connected.
Let $V$ and $W$ be two objects of $\ms U$ over $k$.  A general map
$V\to W(n)$ has as cokernel $Q$ an invertible twisted sheaf supported
on a finite \'etale subscheme of $\ms C^\circ$.  The uniformity
conditions ensure that $Q$ moves in an irreducible family $I$.  The
vector bundle over $I$ parametrizing extensions of $Q$ by $V$ then
surjects onto $\ms U$, showing that it is connected.  Of course, this
outline is hopelessly na\"ive, but there is a version of it that 
works, which we now describe.

\begin{defn}\label{sec:moduli-unif-twist-6}
  Let $\mc C\to\ms G$ be a finite flat covering by a scheme.  The \emph{$\mc
  C$-regularity\/} of a $\ms G$-twisted sheaf $\ms F$ is the
Castelnuovo-Mumford regularity of $\ms F|_{\mc C}$ with respect to the
polarization pulled back from $C$.
\end{defn}

\begin{lem}\label{sec:moduli-unif-twist-7}
  Given a positive integer $N$, there is a quasi-compact open substack
  $\ms U_N\subset\ms U$ parametrizing uniform locally free $\ms
  G$-twisted sheaves of $\mc C$-regularity at most $N$.  Moreover, the
  natural inclusion $\cup_N\ms U_N\inj\ms U$ is an isomorphism.
\end{lem}
\begin{proof}
  By Proposition 3.2.2 of \cite{orbimoduli}, the pullback morphism $\ms
  U\to\Sh(\mc C)$ is of finite presentation, where $\Sh(\mc C)$
  denotes the stack of flat families of locally free sheaves on $\mc
  C$.  On the other hand, we know that the open substack of $\Sh(\mc
  C)$ parametrizing sheaves of regularity at most $N$ is
  quasi-compact.  We conclude that $\ms U_N$ is quasi-compact, as desired.
\end{proof}

Since $\ms U$ is the union of the quasi-compact open substacks $\ms U_N$,
it suffices to prove that each $\ms U_N$ is connected.  Thus, \emph{in what
follows we will fix $N$ and write $\ms U$ in place of $\ms U_N$\/}.

\begin{prop}\label{sec:moduli-unif-twist-8}
  There exists a positive integer $n$ such that for any two objects
  $V$ and $W$ in $\ms U_{spec k}$, a general map $V\to W(n)$ has cokernel $Q$
  satisfying the following conditions.
  \begin{enumerate}
  \item The support $S:=\supp Q$ is a finite reduced substack of $\ms
    G^\circ$ and $Q$ is identified with a $\ms G_S$-twisted invertible
    sheaf.
  \item For any $i$ we have $|S\cap C_i|=\ell^2nH\cdot C_i$.
  \item For $i=1,\ldots,s$, there is a partition $S\cap
    C_i=S_0\coprod\cdots\coprod S_{\ell-1}$ such that (1) for each
    $j$, $|S_j|=\ell n H\cdot C_i$, and (2)
    $Q|_{S_j}\cong\Lambda_{S_j}\tensor\chi^j$.
  \end{enumerate}
\end{prop}
In other words, the cokernel of a general map of uniform sheaves is
itself ``uniformly distributed'' with respect to the characters of
$\m_\ell$.
\begin{proof}
  Using methods similar to those of Lemma 3.1.4.8ff of
  \cite{twisted-moduli}, for any $n$ we have an exact
  sequence of sets
\begin{equation}\label{E:diag}
  \xymatrix{\Hom(V,W(n))\ar[r]&\prod_i\Hom(V|_{C_i},W(n)|_{C_i})\ar@<1ex>[r]\ar@<-1ex>[r]&\prod_j\Hom(V|_{\xi_j},W(n)|_{\xi_j}).}
\end{equation}
For $n$ sufficiently large, we have that the composition
$\Hom(V,W(n))\to\prod\Hom(V|_{\xi},W(n)|_{\xi})$ and the projections
$\Hom(V,W(n))\to\Hom(V|_{C_i},W(n)|_{C_i})$ are surjective.  Fixing
identifications $V|_{\xi_j}\simto V(j)$ and $W|_{\xi_j}\simto V(j)$,
we may consider only those maps $V\to W(n)$ and $V|_{C_i}\to
W(n)|_{C_i}$ which induce the identity on $V(j)$ for each $j$; these
form linear subspaces which we will denote $H$ and $H_i$,
respectively.  For sufficiently large $n$, the natural maps $H\to H_i$
are surjective.

We first claim that it suffices to prove the statement for the
restrictions $V_i$ and $W(n)_i$ to $C_i$.  If this statement is proven
for each $i$ then we can use Diagram \eqref{E:diag} to glue such maps
together.  Since the conditions are clearly open, this will show that
a general map $V\to W(n)$ has this property.  Thus, we will restrict
our attention to $C_i$ in what follows.

If $i=1,\ldots,s$, then $V$ and $W$ have eigendecompositions, and maps
between $V$ and $W(n)$ must preserve the eigensheaves, each of which
has rank $\ell$.  Thus, to show the result we may assume that for
$i=1,\ldots,s$ the sheaves $V$ and $W$ are generically isomorphic
eigensheaves of degree $0$ and rank $\ell$.  Suppose that
$\Lambda_i^{\vee}\tensor W$ is isotyptic of type $\chi^j$.  Then any
simple quotient of $W(n)$ supported on a reduced residual gerbe $x$ in
$\ms C_i^\circ$ will be of the form $(\Lambda_i)_x\tensor\chi^j$, and
we see that to prove the result it suffices to prove that a general
map between $V$ and $W(n)$ which induces an isomorphism at each
$\xi_j$ has cokernel which is an invertible twisted sheaf supported on
a finite reduced subscheme of $C_i^\circ$.  

Twisting down by $\Lambda_i\tensor\chi^j$, we reduce to proving the
following statement.  Let $C$ be a smooth projective curve over an
algebraically closed field and let $D\subset C$ be a reduced effective
Cartier divisor.  Given two locally free sheaves $V$ and $W$ of rank
$\ell$ on the root stack $C\{D^{1/\ell}\}$ such that $\det V\cong\det
W\cong\ms O_C$ and $V_D\cong W_D$, a general map $V\to W(n)$ with
prescribed value over $D$ has cokernel an invertible sheaf supported
on a finite reduced subscheme of $C\setminus D$.  

Let $\Phi:V|_H\to W(n)_H$ be the universal map on $C\times H$.  We are
interested in the locus over $H$ over which $\coker\Phi$ is smooth
over $C$ with fibers of rank $\leq 1$.  This is clearly an open
condition on $H$, and the complement of this locus forms a cone $B$.  We
will investigate the fiber dimension of this cone over $C$.

Given a point $c\in C^\circ$, we have that the restriction map
$$H\to\hom(V_{\spec(\ms O_{C/c}/\mf m_c^2)},W(n)_{\spec(\ms O_{C/c}/\mf
m_c^2)})=\M_\ell(\ms O_{C,c}/\mf m_c^2)$$ is surjective.  Writing a matrix as $A+\eps B$, we have by the Jacobi formula that 
$$\det (A+\eps B)= \det A+\eps\Tr(\operatorname{adj}(A)B),$$
where $\operatorname{adj}$ denotes the classical adjoint (arising from
cofactor matrices).  The vanishing of this expression is clearly of
codimension at least $3$, so the fiber of $B$ over $c$ is of
codimension at least $3$ in $H$.  This shows that $B$ has codimension
at least $2$ in $H$, whence $H$ contains a point $h$ such that
$\coker\Phi_h$ is smooth.  To show that a general such point has fiber
rank $1$, it suffices by a similar argument to show that the locus of
element in $\M_\ell(k)$ in which the rank drops by at least $2$ has
codimension at least $2$.  This can be settled by a trick: a
perturbation argument shows that the locus in question is the boundary
of the locus of non-injective maps, which is itself of codimension
$1$.
\end{proof}

For each $i$, define a $k$-stack $\ms Q_i$ as follows.  The objects of
$\ms Q_i$ over $T$ are pairs $(E,L)$ with $E\subset C_i^\circ\times T$
a closed subscheme which is finite \'etale over $T$ of degree $n^2\ell
H\cdot C_i$ and $L$ an invertible $\ms G_E$-twisted sheaf such that if
$i=1,\ldots,s$, then there is a partition $E=E_0\coprod\cdots\coprod
E_{\ell-1}$ with $L|_{E_j}\tensor\Lambda_i^{\vee}$ isotypic of type
$\chi^j$.  (If $i>s$ there is no condition, and $\ms Q_i$ parametrizes
pairs $(E,L)$ with $L$ any invertible $\ms G_E$-twisted sheaf.)

\begin{prop}\label{sec:moduli-unif-twist-9}
  The stack $\ms Q_i$ is irreducible.
\end{prop}
\begin{proof}
  First suppose $i=1,\ldots,s$, so that $R_i\to C_i^\circ$ is split.
  We claim that $\ms Q_i$ is a torus-gerbe over the dense open
  subscheme of $\prod_{j=0}^{\ell-1}\left(\Sym^{n\ell H\cdot
      C_i}C_i^\circ\right)$ parametrizing (unordered) tuples of
  pairwise distinct points.  Indeed, given a pair $(E,L)$ with
  partition $E=E_0\coprod\cdots\coprod E_{\ell-1}$, the partition
  defines a point of the product of symmetric powers which lies in the
  locus parametrizing pairwise distinct points.  Moreover, this point
  completely characterizes the partition of $E$ (including the weight
  assigned to the restriction of $L$ to each element of the
  partition).  Thus, the fibers of this map are precisely given by
  stacks of line bundles on each $E_i$.  These are gerbes banded by
  tori of rank $n\ell H\cdot C_i$, since each $E_i$ consists of $\ell$
  points.  Since the symmetric powers are irreducible, it follows that
  $\ms Q_i$ is irreducible.

  Now suppose that $i>s$.  We claim that $\ms Q_i$ is a gerbe over a
  dense open subscheme of $\Sym^{n^2\ell H\cdot C_i}R_i$ banded by a
  torus of rank $n^2\ell H\cdot C_i$.  To see that there is a map $\ms
  Q_i\to\Sym^{n^2\ell H\cdot C_i}R_i$, note first that if
  $E=T\coprod\cdots\coprod T\inj (C_i^\circ)_T$ is the trivial finite \'etale
  cover then the invertible sheaf $L$ yields $n^2\ell H\cdot C_i$ maps
  from $T$ to $R_i$, as $R_i$ is the relative twisted Picard space.
  By composition, this gives a point of $\Sym^{n^2\ell H\cdot
    C_i}R_i$.  Moreover, acting by an automorphism of the pair
  $(T\coprod\cdots\coprod T\inj (C_i^\circ)_T,L)$ (i.e., an automorphism of $L$)
  does not change the image point in $\Sym^{n^2\ell H\cdot C_i}R_i$.
  By descent theory, any pair $(E,L)$ gives rise to a point of
  $\Sym^{n^2\ell H\cdot C_i}R_i$, as desired.  Since the (geometric)
  points of $E$ are pairwise distinct in $C_i^\circ$, the image of
  $(E,L)$ lies in the open subscheme $O$ of $\Sym^{n^2\ell H\cdot C_i}R_i$
  parametrizing points with pairwise distinct images in $C_i$.

  The statement that $\ms Q_i$ is a torus gerbe over $O$ can be
  checked \'etale-locally on $T$, so we may assume that
  $E=T\coprod\cdots\coprod T$.  Suppose two pairs
  $(T\coprod\cdots\coprod T\inj (C_i^\circ)_T,L_1)$ and
  $(T\coprod\cdots\coprod T\inj (C_i^\circ)_T,L_2)$ give rise to the
  same point of $\Sym^{n^2\ell H\cdot C_i}R_i$.  Since each summand
  $T$ maps to a distinct section of $(C_i^\circ)_T$, it follows that
  upon possibly reordering the summands, for each summand $T$ the two
  restrictions $(L_1)_T$ and $(L_2)_T$ must map to the same section of
  $R_i$ and thus differ by an invertible sheaf on $T$.  Shrinking $T$
  if necessary, we see that any two pairs with the same image in
  $\Sym^{n^2\ell H\cdot C_i}R_i$ are locally isomorphic.  In addition,
  since any point of $\Sym^{n^2\ell H\cdot C_i}R_i$ parametrizing
  distinct points comes \'etale locally from a finite collection of
  disjoint sections of $R_i$ and $R_i$ is the relative Picard space,
  we see that any section of $O$ locally gives a point of $\ms Q_i$.
  Thus, $\ms Q_i\to O$ is a gerbe.  Finally, the band of $\ms Q_i$ is
  a torus because automorphisms of a pair $(E,L)$ are given by
  automorphisms of $L$, which is an invertible sheaf over a finite
  \'etale $T$-scheme of rank $n^2\ell H\cdot C_i$.

  Since $R_i$ is connected and smooth, it is irreducible, whence its
  symmetric powers are irreducible.  But then, since a torus $T$ is
  irreducible, we see that any $T$-gerbe over $\Sym^{n^2\ell H\cdot
    C_i}R_i$ is also irreducible, as desired.
\end{proof}

Let $\ms Q=\prod\ms Q_i$.  Define a geometric vector bundle over $\ms
Q$ as follows.  Let $\mf Q$ be the universal object on $\ms G\times\ms
Q$.

\begin{lem}\label{sec:moduli-unif-twist-10}
  The complex $\R(\pr_2)_\ast\rshom(\mf Q,\pr_1^\ast V)[1]$ is
  quasi-isomorphic to a locally free sheaf $\ms F$ on $\ms Q$.
  Moreover, his sheaf has the property that for any affine scheme $T$
  and any morphism $\psi:T\to\ms Q$, the set $\ms F_T(T)$ parametrizes
  extensions $0\to V\to W\to Q\to 0$ with $Q$ the object of $\ms Q$
  corresponding to $\psi$.
\end{lem}
\begin{proof}
  By standard nonsense, the complex $\R(\pr_2)_\ast\rshom(\mf
  Q,\pr_1^\ast V)$ is perfect and commutes with arbitrary base change
  on $\ms Q$.  Moreover, since $\mf Q$ has torsion fibers (over
  $\pr_2$) and $V$ is locally free, we have that for each geometric
  point $q:\spec\widebar\kappa\to\ms Q$ the sheaf $h^0(\L
  q^\ast\R(\pr_2)_\ast\rshom(\mf Q,\pr_1^\ast V))$ is zero, and
  similarly for $h^i$ with $i>1$.  It follows that
  $\R(\pr_2)_\ast\rshom(\mf Q,\pr_1^\ast V)$ is a locally free sheaf
  concentrated in degree $1$ of formation compatible with arbitrary
  base change on $\ms Q$, with affine fibers the groups  
  $\ext^1(Q,V)$, as claimed.
\end{proof}

Thus, we see that locally free extensions of objects of $\ms Q$ by $V$ are
parametrized by a dense open substack of the irreducible stack $\V(\ms
F^\vee)$.  Let $U\subset\V(\ms F^\vee)$ be the open stack parametrizing sheaves
of the form $W(n)$, where $W$ is a uniform
$\ms G$-twisted sheaf.  

\begin{proof}[Proof of Theorem \ref{sec:moduli-unif-twist-1}]
  First suppose $\ms U$ is non-empty.  By Lemma
  \ref{sec:moduli-unif-twist-8}, $U$ is non-empty and the map $U\to\ms
  U$ is surjective.  Since $U$ is irreducible, we conclude that
  $\ms U$ is irreducible, as desired.

  Thus, it suffices to prove that $\ms U$ is non-empty, i.e., that
  there exists a uniform $\ms G$-twisted sheaf of rank $\ell^2$ and
  determinant $\ms O$.  This is somewhat subtle, and is carried out in
  Paragraph \ref{sec:moduli-unif-twist-11} below.
\end{proof}

\begin{para}\label{sec:moduli-unif-twist-11}
  In this paragraph, we will show that $\ms U$ is non-empty.  Oddly
  enough, it is here that we must assume $\ell$ is odd.  Fix totally biregular
  $\ms G(j)$-twisted sheaves $V(j)$ for each singular gerbe $\xi_j$.  
\end{para}

\begin{lem}\label{L:one-at-a-time}
  It suffices to prove that for each $i$ there is a uniform locally
  free $\ms G_i$-twisted sheaf $V$ of rank $\ell^2$ with trivial
  determinant such that for each singular gerbe $\xi_j$ contained in
  $\ms C_i$, we have that $V|_{\xi_j}\cong V(j)$.
\end{lem}
\begin{proof}
  Indeed, we can glue the sheaves $V$ together at the nodes of $\ms C$
  to get a uniform sheaf on the entire residual curve.  (This is
  trivial because $k$ is now assumed to be algebraically closed.)
\end{proof}

\begin{prop}
  The stack $\ms U$ is nonempty.
\end{prop}
\begin{proof}
  Using Lemma \ref{L:one-at-a-time}, we may assume $C=C_i$ for a fixed
  $i$. If the ramification cover $R_i\to C_i$ is connected then by
  definition we just need a locally free $\ms G_i$-twisted sheaf with
  the given local structures $V(j)$ and trivial determinant.  Choosing
  any locally free $\ms G_i$-twisted sheaf $W$ of rank $\ell^2$ with
  totally regular generic fiber, we can glue on the local structures
  $V(j)$ as follows.  Write $x_j\in C_i$ for the image of $\xi_j$.
  Since $\m_\ell$ is reductive, infinitesimal deformation of $V(j)$
  over $\spf\widehat{\ms O}_{C_i,x_j}$ is unobstructed.  Since $\ms
  G_i\to C_i$ is proper, we can apply the Grothendieck existence
  theorem to any such infinitesimal obstruction to yield a locally
  free $\ms G_i\tensor\spec(\widehat{\ms O}_{C_i,x_j})$-sheaf
  $\widehat V$ whose restriction to $\xi_j$ is isomorphic to $V(j)$.
  On the other hand, writing $K$ for the fraction field of
  $\widehat{\ms O}_{C_i,x_j}$, we have that $W_K$ and $\widehat V_K$
  are isomorphic (as there is only one totally regular structure of a
  given rank over a field).  Choosing a generic isomorphism and
  applying Corollary \ref{C:formal extension} yields a locally free
  $\ms G_i$-twisted sheaf $V$ which is totally biregular and isomorphic
  to $V(j)$ over $\xi_j$ for each $j$.

  It remains to handle the case in which $R_i\to C_i$ is split.  In
  this case, since $\ms G$ is geometrically trivial at residual gerbes
  (by assumption), there is an invertible $\ms G_i$-twisted sheaf
  $\Lambda$ such that $\Lambda^{\tensor\ell}$ is the pullback of an
  invertible sheaf $\lambda$ on $C_i$.  Let
  $V_0,V_1,\ldots,V_{\ell-1}$ be locally free sheaves on $\ms C_i$
  such that
  \begin{enumerate}
  \item for each $s=0,\ldots,\ell-1$ the sheaf $V_s$ is generically
    isotypic on the $j$th power of the generic stabilizer $\m_\ell$ of
    $\ms C_i$,
  \item for each $j$ the restriction $(V_s)(j)$ is isomorphic to the
    representation of $\m_\ell\times\m_\ell$ which is the regular
    representation of the second factor tensored with the $s$-power
    character of the first factor, and
  \item the determinant of $\Lambda\tensor V_s$ is trivial for each $s$.
  \end{enumerate}
  This is possible because $\ell$ is odd, so the regular
  representation has trivial determinant.  The sheaf
  $\Lambda\tensor(V_0\oplus\cdots\oplus V_{\ell-1})$ is then the
  desired uniform $\ms G_i$-twisted sheaf of rank $\ell^2$.
\end{proof}

\section{Curves over higher local fields}
\label{S:local-mama}

In this section we fix a positive integer $M$.

Let $R$ be an excellent Henselian discrete valuation ring with fraction field
$K$, uniformizer $t$, and residue field $k$ of characteristic exponent
$p$.  Suppose $n$ is invertible in $k$ and that $R$ contains a
primitive $n$th root of unity.

\begin{lem}\label{perfection lemma} There is an inclusion $R\subset
  R'$ such that
  \begin{enumerate}
  \item $R'$ is an excellent Henselian discrete valuation ring with uniformizer
    $t$;
  \item $R'$ is a colimit of finite free extensions of $R$ of
    $p$-power degree;
  \item the morphism of residue fields $k\inj k'$ identifies $k'$ with
    a perfect closure of $k$.
  \end{enumerate}
\end{lem}
\begin{proof} If $p=1$, we $R'=R$; for the rest of the proof, we
  assume that $p>1$.  Given an element $\bar x\in k\setminus k^p$, we
  can choose a lift $x\in R$ and consider the extension
  $S=R[y]/(y^p-x)$.  It is easy to see that $S$ is finite and free
  over $R$ of degree $p$ and that $t$ generates a maximal ideal, so
  that $S$ is a dvr with uniformizer $t$.  The result follows by
  transfinite induction.
\end{proof}

\begin{lem}\label{L:nodal-curve}
  If $k$ has $M'$-Brauer dimension at most $a$ and $C$ is a nodal
  (possibly non-proper) curve over $k$, then for any $\m_n$-gerbe $\ms
  S\to C$ with $n\in M'$, there is a locally free $\ms S$-twisted sheaf of rank $n^a$.
\end{lem}
\begin{proof}
  The conclusion holds when $C$ is smooth by the definition of the
  Brauer dimension and the fact that a torsion-free coherent
  twisted sheaf is locally free (so that any coherent torsion-free
  extension of a generic twisted sheaf will satisfy the conclusion of
  the lemma).  Thus, the conclusion holds for the normalization of
  $C$.  Moreover, at a node of $C$, the restrictions of the locally
  free $\ms S$-twisted sheaves on the two branches are isomorphic.  It
  follows by a basic descent argument that the locally free sheaves on
  the components of $C$ may be glued at the nodes to achieve the
  desired result.
\end{proof}

\begin{thm}\label{T:main} 
  If $k$ has $M'$-Brauer dimension at most $a$ then $K$ has $(Mp)'$-Brauer
  dimension at most $a+1$.
\end{thm}
\begin{proof} We have to check two things: (1) for any finite extension
  $J$ of $K$ and any $\alpha\in\Br'(J)$ we have that
  $\ind(\alpha)|\per(\alpha)^a$, and (2) for any finitely generated
  extension $J(D)$ of transcendence degree $1$ and any
  $\alpha\in\Br'(J(D))$ we have $\ind(\alpha)|\per(\alpha)^{a+1}$.

  First, let $J$ be a finite extension of $K$ with ring of integers
  $S$ and uniformizer $u$.  Given a Brauer class $\alpha\in\Br(J)$ of
  period $n$ prime to $p$, it is clear that $\alpha$ will be unramified
  over $S[u^{1/n}]$, which is finite free of degree $n$ over $S$.  Let
  $\ms X\to\spec S[u^{1/n}]$ be a $\m_n$-gerbe representing $\alpha$.
  The residue field $L$ has Brauer dimension at most $a$ by hypothesis,
  so that there is a locally free $\ms X\tensor L$-twisted sheaf of
  rank $n^{a-1}$.  Deformation theory lifts this to a locally free
  $\ms X$-twisted sheaf of rank $n^{a-1}$, which pushes forward at the
  generic point to yield a locally free twisted sheaf of rank $n^{a}$
  on a gerbe over $\spec J$ representing $\alpha$.  This completes the
  proof of (1).

  It remains to verify (2).  Any finitely generated extension of
  transcendence degree $1$ has the form $J(D)$ with $J$ finite over
  $K$ and $D$ a geometrically integral curve over $J$.  Given
  $\alpha\in\Br'(J(D))$, we may normalize $R$ in $J$ and thus assume
  that $J=K$.  We may also assume that $\per(\alpha)=n$ is prime (and
  relatively prime to $p$ by assumption).  By the standard
  finite-presentation tricks, we see that we can replace $R$ with the
  $R'$ of Lemma \ref{perfection lemma} and thus assume that $k$ is perfect.

  Let $X\to\spec R$ be a regular model of $D$ over $R$, which exists
  by \cite{lipman}.  After blowing up $X$ if necessary, we may assume
  that the union of the ramification divisor of $\alpha$ and the
  (reduced) special fiber of $X/R$ is a strict normal
  crossings divisor.  Applying Corollary \ref{C:unram-stack-I}, we see that
  there is an orbisurface $\ms X/T$ with a generically finite map $\ms
  X\to X$ of degree $n$ over which $\alpha$ is unramified.  It thus
  suffices to show that the index of $\alpha$ restricted to $\ms X$
  divides $n^{a}$.   Let $\ms G\to\ms X$ be a
  $\m_n$-gerbe representing $\alpha$ and write $\ms X_0$ for the
  reduced structure on the special fiber of $\ms X$ over $T$.  

  In order to prove (2), it suffices to find a locally free $\ms
  G$-twisted sheaf of rank $n^{a}$.  Since $\ms X_0$ is
  $1$-dimensional and tame, we see from Proposition \ref{P:deformation theory}
  that any locally free $\ms G|_{\ms X_0}$ twisted sheaf of rank $n^a$
  has vanishing obstruction space and thus deforms to a locally free
  formal $\ms G$-twisted sheaf.  Applying the Grothendieck existence
  theorem (Proposition \ref{P:deformation theory}(3)), we can effectivize the
  deformation to yield a locally free $\ms G$-twisted sheaf of rank
  $n^a$.  Thus, it suffices to show that there is a locally free $\ms
  G|_{\ms X_0}$-twisted sheaf of rank $n^a$.

  Write $\pi:\ms X_0\to Y$ for the coarse moduli space of the stack
  $\ms X_0$.  There are finitely many closed points $y_1,\ldots,y_r\in
  Y$ such that (1) $\pi$ is an isomorphism over
  $Y\setminus\{y_1,\ldots,y_r\}$ and (2) the reduced structure
  $\xi_i:=(\ms X_0\times_Y y_i)_{\text{red}}$ is isomorphic to
  $\B{\m_{n,\kappa(y_i)}}$ for each $i$.  Applying
  Corollary \ref{C:stack-pt-det} and the hypothesis that $k$ has Brauer
  dimension at most $a$, we see that for each $i=1,\ldots,r$, there is
  a locally free $\ms G|_{\xi_i}$-twisted sheaf $\widebar{\ms V}_i$ of
  rank $n^a$.  Moreover, since $\m_n$ is reductive, it follows from
  Proposition \ref{P:deformation theory} that infinitesimal deformations of
  $\widebar{\ms V}_i$ are universally unobstructed (the obstruction
  space vanishes).  Thus, writing $\ms X_0^{(i)}=\ms X_0\times_Y\spec
  \widehat{\ms O}_{Y,y_i}$, it follows from Proposition \ref{P:deformation
    theory}(3) that for each $i=1,\ldots,r$, there is a locally free
  $\ms G_{\ms X_0^{(i)}}$-twisted sheaf of rank $n^a$.

  On the other hand, $C:=\ms X_0\setminus\{\xi_1,\ldots,\xi_s\}\cong
  Y\setminus\{y_1,\ldots,y_r\}$ is a nodal curve over $k$, so by
  Lemma \ref{L:nodal-curve} there is a locally free $\ms G_C$-twisted sheaf of
  rank $n^a$.  Applying Corollary \ref{C:formal extension}, we see that
  we can glue the twisted sheaf over $C$ to the twisted sheaves
  on each $\ms X_0^{(i)}$ to find a locally free $\ms G|_{\ms X_0}$-twisted
  sheaf of rank $n^a$, as desired.
\end{proof}

\begin{cor} Let $C$ be a regular curve over a field $K$ and
  $\alpha\in\Br(K(C))$.
  \begin{enumerate}
  \item If $K$ is a $d$-local field with finite residue field of
    characteristic $p$ and $\per(\alpha)\in p'$ then
    $\ind(\alpha)|\per(\alpha)^{d+1}$.
  \item If $K$ is a $d-1$-local with residue field of transcendence
    degree $1$ over an algebraically closed field (e.g.,
    $K=k(x_1)\(x_{2}\)\cdots\(x_d\)$) of characteristic exponent $p$
    and $\per(\alpha)\in p'$, then $\ind(\alpha)|\per(\alpha)^d$.
  \item If $K$ is an algebraic extension of a local field of residual
    characteristic exponent $p$ which contains the maximal unramified
    extension and $\per(\alpha)\in p'$, then
    $\ind(\alpha)=\per(\alpha)$.
  \end{enumerate}
\end{cor}
\begin{proof} The first two statements follow by induction, starting
  with class field theory or de Jong's theorem \cite{dejong-per-ind},
  while the third is a direct application of the theorem (after
  reducing to the case where $K$ is finite over the maximal unramified
  extension).
\end{proof}

\section{Curves over global fields}
\label{sec:curves-over-global-1}

Let $C/K$ be a proper smooth curve over a global field with a point
$P\in C(K)$.  In this section we link the period-index problem for
certain classes in $\Br(C)$ to the Hasse principle for geometrically
rational varieties over $K$.

We first recall a well-known conjecture of Colliot-Th\'el\`ene.  What
is written here is a special case.  For the general case (which just
says that the Brauer-Manin obstruction is the only one for $0$-cycles
of degree $1$ on smooth proper varieties) the reader is referred to Conjecture 2.5 of \cite{MR1743234}.

\begin{conj}[Colliot-Th\'el\`ene]\label{sec:curves-over-global-5}
  Suppose $V/K$ is a smooth projective variety such that
  $V\tensor\widebar K$ is rational and $\Pic(V\tensor\widebar K)=\Z$.
  If there is a $0$-cycle of degree $1$ on $V\tensor K_{\nu}$ for all
  places $\nu$ of $K$ then there is a $0$-cycle of degree $1$ on $V$.
\end{conj}

We will link this to the period-index problem as follows.  Write
$\Br_P(C)$ for the kernel of the restriction map
$P^\ast:\Br(C)\to\Br(K)$.  Given an integer $n$, say that a place of
$K(C)$ \emph{divides $n$\/} if the restriction of the valuation on
$K(C)$ to the prime field has residue field of characteristic dividing
$n$.  (In geometric language, the center of the valuation on a proper
model of $C$ over a scheme of integers of $K$ has image contained in
the closed subscheme cut out by the function $n\cdot 1$.)

\begin{thm}\label{sec:curves-over-global-6}
  If Conjecture \ref{sec:curves-over-global-5} holds, then any
  $\alpha\in\Br_P(C)$ of odd period which is unramified at all places
  dividing $\per(\alpha)$ satisfies $\ind(\alpha) | \per(\alpha)^2$.
\end{thm}
\begin{proof}[Proof of Theorem \ref{sec:curves-over-global-6} when $\alpha$ has prime period $\ell$]
  Let $\per(\alpha)=\ell$.  If $K$ has characteristic $0$, let $T$ be
  the scheme of integers in $K$; if $K$ has positive characteristic,
  choose a transcendence basis for $K$ and let $T$ be the
  normalization of $\P^1$ in the resulting function field extension.
  We will call $T$ ``the'' scheme of integers of $K$.  Note that when
  $K$ has infinite places $\infty$, the assumption that $\per\alpha$
  is odd implies that $\alpha|_{C_\infty}=0$.

  Taking a branched cover if necessary, we may assume that $C$ has
  genus at least $2$.  Choose a regular proper model $\mc C\to T$ and
  let $D\subset\mc C$ be the ramification divisor of $\alpha$.  Since
  $\alpha\in\Br(C)$, we see that $D$ is supported in finitely many
  fibers of $\mc C\to T$, and the assumption on the ramification not
  dividing the period of $\alpha$ shows that the residue fields of the
  components of $D$ have characteristics prime to $\ell$.  Using
  Lipman resolution \cite{lipman}, we may assume that the fibers containing
  the ramification divisor of $\alpha$ form a snc divisor which we
  will abusively write as $D\subset\mc C$.  Applying
  Proposition \ref{sec:splitt-root-constr-1}, we see that $\alpha$ extends to a
  class in $\Br(C\{D^{1/\ell}\})$.  Moreover, the scheme $D$ is a
  union of finitely many residual curves $D_1,\ldots,D_b$.  Choosing a
  $\m_\ell$-gerbe $\ms G\to\mc C\{D^{1/\ell}\}$ such that $\ms
  G|_{C_{\widebar K}}=0$, we see that for each $i=1,\ldots,b$, there
  is a locally free $\ms G|_{D_i}$-twisted sheaf of trivial
  determinant and rank $N\ell^2$ with $\ell\not|N$.  Taking direct
  sums if necessary, we may assume that $N$ is the same for all $i$.

  Let $\ms V\to\spec K$ be the stack of stable locally free $\ms
  G_K$-twisted sheaves of rank $N\ell^2$ with determinant $\ms O(P)$.
  It follows from Proposition 3.1.1.6 of \cite{twisted-moduli} and
  \cite{king-schofield} that $\ms V$ is a $\m_{N\ell^2}$-gerbe over a
  smooth projective geometrically rational $K$-variety $V$ with
  $\Pic(V\tensor\widebar K)=\Z$.  We claim that there exists an $N$
  relatively prime to $\ell$ such that for all places $\nu$ of $K$, we
  have $\ms V(K_{\nu})\neq \emptyset$.  Assuming we have done this, we
  can then apply Conjecture \ref{sec:curves-over-global-5} to conclude
  that $V$ has a $0$-cycle of degree $1$.  Moreover, since we will
  prove that the stack (and not merely the coarse space) has local
  points, it follows from a simple argument with the Leray spectral
  sequence that $\Br(K)\to\Br(V)$ is surjective and then that the
  Brauer class of $\ms V\to V$ is trivial.  Thus, the $0$-cycle of
  degree $1$ on $V$ lifts to $\ms V$.  But then there is a locally
  free $\ms G_K$-twisted sheaf of rank $M\ell^2$ for some $M$
  relatively prime to $\ell$, so that $\ind(\alpha) | M\ell^2$.  Since
  $\per(\alpha)=\ell$, we conclude by standard methods that
  $\ind(\alpha)|\ell^2$, as desired.

  Thus, it remains to find local objects in $\ms V(K_\nu)$.  Since
  adjoining a primitive $\ell$th root of unity is a prime-to-$\ell$
  extension, we may assume by Lemma \ref{sec:moduli-unif-twist-3} that
  $K$ contains a primitive $\ell$th root of unity.  If $\nu$ is not in
  the image of the ramification locus (including infinite $\nu$), we
  claim that $\alpha|_{K_\nu}=0$. If $\nu$ is infinite this is an
  elementary parity argument.  If $\nu$ is finite, then this follows
  from the fact that the Brauer group of a proper curve over a
  complete dvr with finite residue field is finite.  In these cases,
  the existence of an object of $\ms V(K_\nu)$ follows
  immediately from the fact that stable vector bundles of arbitrary
  invariants exist on curves of genus at least $2$ over any infinite
  field.

  Now suppose $\nu$ has center under one of the residual curves $D_i$.
  By Corollary \ref{sec:moduli-unif-twist-2}, there is some $N$
  relatively prime to $\ell$ and a locally free $\ms G|_{D_i}$-twisted
  sheaf $V_i$ of rank $\ell^2$ and trivial determinant.  Taking
  suitable direct sums, we may assume that $N$ is independent of $i$.
  The deformation theory of twisted sheaves in Proposition
  \ref{P:deformation theory} shows that $V_i$ deforms to yield a
  locally free $\ms G\tensor K_\nu$-twisted sheaf $\widetilde V_i$ of
  rank $N\ell^2$ with trivial determinant.  Since $\alpha|_P=0$, we
  can make an elementary transform to yield a locally free $\ms
  G\tensor K_\nu$-twisted sheaf $W_i$ of rank $N\ell^2$ and
  determinant $\ms O(P)$.  Let $\ms S$ denote the stack of all locally
  free $\ms G\tensor K_\nu$-twisted sheaves of rank $N\ell^2$ and
  determinant $\ms O(P)$.  We know that $\ms V\subset\ms S$ is a dense
  open immersion of irreducible smooth stacks.  Thus, any algebraic
  deformation of $W_i$ over a smooth base which is versal at $[W_i]$
  will contain a dense open parametrizing stable sheaves.  But any
  smooth connected $K_\nu$-variety with a rational variety has a dense
  set of rational points, so we conclude that there is a rational
  point parametrizing a stable twisted sheaf, yielding an object of
  $\ms V(K_\nu)$, as desired.
\end{proof}

\begin{proof}[Proof of Theorem \ref{sec:curves-over-global-6} in the general case]
  The Schur decomposition of $\alpha$ immediately reduces to the
  situation where $\per(\alpha)$ is a prime power, say
  $\per(\alpha)=\ell^n$.  The class $\alpha':=\ell^{n-1}\alpha$ has
  period $\ell$, so assuming that we have treated the case of prime
  period (done below), we can find a covering $X'\to X$ of degree
  $\ell^2$ splitting $\alpha'$.  Since the pullback of $\alpha$ to
  $X'$ has period $\ell^{n-1}$, we are done by induction, given the
  following (crucial) lemma.
  \begin{lem}
    If $\beta\in\Br(\ms O_{X,P})$ is a class whose restriction to $P$
    is trivial, then there is a splitting field $K'/K(X)$ such that
    the normalization of $X$ in $K'$ has a $K$-rational point lying over $P$.
  \end{lem}
  \begin{proof}
    Let $\P\to\spec\ms O_{X,P}$ be a Brauer-Severi scheme representing
    $\beta$.  Since $P\tensor\kappa(P)$ is trivial, there is a
    $\kappa(P)$-rational point $Q\in \P_P$.  Let $U$ be an affine
    neighborhood of $Q$ in $\P$.  Choose a regular sequence $\widebar
    f_1,\ldots,\widebar f_d$ cutting out $Q$ in the fiber $\P_P$ and
    choose arbitrary lifts $f_1,\ldots,f_d\in\ms O_{\P}(U)$.  Let $Z$
    be the subscheme of $U$ determined by $f_1,\ldots,f_d$.  By
    construction, $Z$ is flat over $\ms O_{X,P}$ and normal at $Q$
    with residue field $K=\kappa(P)$.  It follows that the generic
    fiber of $Z$ gives rise to the desired splitting field.
  \end{proof}
  Applying the lemma in our situation, we see that there is some
  $X'\to X$ splitting $\alpha'$, possibly of large degree, such that
  there is a rational point $Q\in X'(K)$ over $P$.  This property
  clearly holds for any intermediate curve $X'\to X''\to X$.  Since
  $\alpha'$ has index dividing $\ell^2$, there will be such an
  intermediate curve $X''$ of degree $\ell^2$ over $X$, and we know
  that there is a rational point $Q\in X''(K)$ such that
  $\alpha|_Q=0$.  Thus, we can proceed by induction.
\end{proof}

We can also get the first ``proof'' of the standard conjecture for a
class of $C_3$-fields.

\begin{cor}\label{sec:curves-over-global-4}
  Suppose $Y$ is a proper smooth surface over a finite field $\F_q$.
  Given $\alpha\in\Br(\F_q(Y))'$, if Conjecture
  \ref{sec:curves-over-global-5} holds then
  $\ind(\alpha)|\per(\alpha)^2$.
\end{cor}
\begin{proof}[Sketch of proof]
  We first record an argument due to de Jong for fibering a birational
  model $Y$ so that the ramification of $\alpha$ ends up in a fiber.
  We can replace $\F_q$ with its maximal prime-to-$q$ extension and
  thus assume that the base field is an infinite algebraic
  extension of its finite prime field.

  We may assume that the ramification divisor of $\alpha$ is a strict
  normal crossings (snc) divisor $D=D_1+\cdots+D_m\subset Y$.  For each
  component $D_i$, there is a finite ramification extension $D_i'\to
  D_i$.  By the Chebotarev density theory, there are closed points
  $d_i\in D_i$ of arbitrarily high degree such that $D_i'$ is totally
  split over $d_i$.  It is an exercise to check (by examining the
  Henselian local structure of $\alpha$ near $d_i$ as in e.g.\
  Proposition 1.2 of \cite{saltman-fix}) that the ramification divisor
  of $\alpha$ in the blow up $Y$ at $d_i$ is the strict transform of
  $D$.  

  Blowing up enough such $d_i$, we may assume that $D$ is contained in
  a snc divisor $E$ whose components have a negative definite
  intersection matrix.  By Theorem 2.9(B) of \cite{artin-contract},
  there is a morphism to a normal projective surface $Y\to \widebar Y$
  which is an isomorphism outside of $E$.  In particular, taking a
  pencil of very ample divisors on $\widebar Y$ which meet
  transversely in the smooth locus, we see that a blow-up of $\widebar
  Y$ in finitely many smooth points fibers over $\P^1$ with a section
  passing entirely through the smooth locus.  Blowing up the
  corresponding locus in $Y$ yields a blow up $\widetilde Y\to Y$ and
  a proper flat morphism $\pi:\widetilde Y\to\P^1$ with a section
  $P:\P^1\to\widetilde Y$ such that
  \begin{enumerate}
  \item $\pi$ has smooth geometrically connected generic fiber;
  \item the ramification of $\alpha$ is entirely contained in a fiber of $\pi$;
  \item $\alpha|_P=0$.
  \end{enumerate}
  The last statement follows from the fact that any base point of the
  pencil gives rise to a section, and that $\Br(\P^1)=0$.  

  Writing $C$ for the generic fiber of $\pi$, we have that
  $\alpha\in\Br_P(C)$.  If $q$ is odd or $q$ is even and the period of
  $\alpha$ is odd then the period-index relation now follows
  immediately from Theorem \ref{sec:curves-over-global-6}.  If $q$ is
  even and $\per(\alpha)$ is $2$ then the the result follows from the
  fact that the absolute Frobenius $F:Y\to Y$ is finite locally free
  of degree $4$ and acts as multiplication by $2$ on the Brauer group.
  The general case follows from the Schur decomposition and induction
  on the power of $2$ dividing the period.
\end{proof}

\appendix

\section{Period-index examples}

\begin{center}
Daniel Krashen\\
Yale University
\vskip 1\baselineskip
\end{center}

In this appendix to the paper of Max Lieblich, we will construct
examples of division algebras over certain fields in order to
demonstrate that the bounds given in corollary \ref{appendix-cor} are
sharp. The main tool in doing this is the use of valuation theory for
division algebras as in \cite{Wad:VT}. I would like to thank A.
Wadsworth for pointing out some results in the literature which helped
to shorten the exposition of these examples.

\subsection{Curves over iterated Laurent series fields}
Let $k$ be an algebraically closed field, $char(k) \neq p$. Let $K =
k\(x_1\)\cdots\(x_r\)(t)$. In this section, we will construct a
division algebra over $K$ with period $p$ and index $p^r$.

Let $K' = k(t)\(x_1\)\cdots\(x_r\)$. Note that $K'$ is a Henselian
valued field with values in $\bZ^r$ and with residue field
$k(t)$. Also note that we have inclusions $k(t) \subset K \subset
K'$. Let $\ell/k(t)$ be abelian Galois with group $(\bZ/p)^r$. By
Kummer theory, we may write $\ell = \ell_1 \otimes_{k(t)} \cdots
\otimes_{k(t)} \ell_r$, with $\ell_i = k(t)(\alpha_i)$, $\alpha_i^p =
a_i \in k(t)$. Using \cite{Wad:VT}, example 3.6, we note that the
inertial lifts $L_i' = K'(\alpha_i)$ of the $\ell_i$ are cyclic
Galois, say with generators $\sigma_i$ of the Galois group, and that
the algebra
$$A' = (L_1'/K', \sigma_1, x_1) \otimes_{K'} \cdots \otimes_{K'}
(L_r'/K', \sigma_r, x_r)$$ is a (nicely semi-ramified) division
algebra. Therefore, if we set $L_i = K(\alpha_i)$ and
$$A = (L_1/K, \sigma_1, x_1) \otimes_{K} \cdots \otimes_{K} (L_r/K,
\sigma_r, x_r),$$ we have $A \subset A \otimes_K K' = A'$. Since $A'$
is division, $A$ contains no zerodivisors and therefore must also be
division. Therefore $A$ is a division algebra of index $p^r$ and
period $p$ as desired.

\subsection{Curves over higher local fields}

In this section we will construct division algebras of period $p$ and
index $p^{r+1}$ over fields of the form $F(t)$, where $F$ is a r-local
field in the following two cases: $F = \bF_q\(x_1\)\cdots\(x_r\)$ with
$q \neq p$ or $F = \bQ_q\(x_1\)\cdots\(x_{r-1}\)$ with $p | (q-1)$.

\begin{lem} \label{basic_lem}
Let $k$ be a global field, and choose $\ell/k$ Galois with group
$(\bZ/p)^r = \left<\sigma_1, \ldots, \sigma_r\right>$. Then we may
choose a cyclic division algebra $D/k$ of index $p$ such that $D
\otimes \ell$ is still division.
\end{lem}
\begin{proof}
To do this we simply note that by Chebotarev density (see
\cite{Serre:ZetaL}, or \cite{Jar}), there are infinitely many primes
such that our field extension $\ell$ is completely split. Let $\fP$ be
such a prime in $k$. By Cassels-Fr\"ohlich page 187 corollary 9.8, a
central simple algebra with a nontrivial Hasse invariant at $\fP$ must
remain nontrivial upon extending scalars to $\ell$. Therefore, we may
simply choose an arbitrary division algebra with a Hasse invariant of
order $p$ at $\fP$. Since the period is equal to the index over a
global field, this algebra also has index $p$.
\end{proof}

Recall that we have either:
\begin{enumerate}
\item \label{case1}
$K = \bF_q\(x_1\)\cdots\(x_r\)(t)$ with $q \neq p$ or 
\item \label{case2}
$K = \bQ_q\(x_1\)\cdots\(x_{r-1}\)(t)$ with $p | (q-1)$.
\end{enumerate}
In the second case, we let $x_r = q$, and in each case $K$ posseses a
valuation into the (lexicographically) ordered abelian group $\bZ^r$
with uniformizers $x_1, \ldots, x_r$. Let $\ov{K} \cong \bF_q(t)$ be
the residue field.

Since $\ov{K}$ is global, we may choose a extension $\ell / \ov{K}$
and a symbol algebra $\cD$ as in lemma \ref{basic_lem}. 

\begin{lem}
We may find an unramified abelian extension $L/K$ of dimension $p^r$
and an inertial symbol algebra $D/K$ of degree $p$ such that $\ov{L} =
\ell$ and $\ov{D} = \cD$.
\end{lem}

Using this lemma, it follows that if we write $L$ as the compositum of
the cyclic algebras $L_i$ with cyclic Galois groups generated by
$\sigma_i$, then the algebra
$$A = D \otimes_K (L_1/K, \sigma_1, x_1) \otimes_K \cdots \otimes_K
(L_r/K, \sigma_1, x_r)$$ is a (valued) division algebra by
\cite{Mor:VD}, theorem 1. Therefore we will have the desired algebra
of index $p^{r+1}$ and period $p$.

\begin{proof}
\begin{case}[1]
$K = \bF_q\(x_1\)\cdots\(x_r\)(t)$.
\end{case}
We have natural inclusions $\ov{K} \subset K \subset
\ov{K}\(x_1\)\cdots\(x_r\)$. It is well known that $\til{L} = \ell \otimes_{\ov{K}}
\ov{K}\(x_1\)\cdots\(x_r\)$ and $\til{D} = \cD \otimes_{\ov{K}}
\ov{K}\(x_1\)\cdots\(x_r\)$ are both division (for example, this
follows from lemma \ref{division_val} below). Since
$\ov{K}\(x_1\)\cdots\(x_r\)$ is Henselian valued, it follows that the valuation
extends to $\til{D}$ and $\til{L}$ and one may check that these are
unramified. Therefore, if we define $L = \ell \otimes_{\ov{K}} K$ and
$D = \cD \otimes_{\ov{K}} K$, we have $L \subset \til{L}$ and $D
\subset \til{D}$, and in particular both are domains and hence
division. These clearly have the desired properties.

\begin{case}[2]
$K = \bQ_q\(x_1\)\cdots\(x_{r-1}\)(t)$ with $p | (q-1)$.
\end{case}
By the assumption on $p$ and $q$, we know that by Hensel's lemma,
$\ov{K}$ and $\bQ_q(t)$ (and hence $K$) contain a primitive $p$'th root
of unity which we will always denote by $\omega$. Therefore Kummer
theory applies and we may write $\ell = \ell_1 \otimes_{\ov{K}} \cdots
\otimes_{\ov{K}} \ell_r$, with each $\ell_i = \ov{K}(\alpha_i)$,
$\alpha_i^p = a_i \in \ov{K}$. Also, we may write $\cD$ as a symbol
algebra $\cD = (a, b)_{\ov{K}, \omega}$.

Choose lifts $\til{a_i}, \til{a}, \til{b} \in \bZ_q(t)^* \subset K$ of
$a_i, a, b$ respectively. Set $L_i = K(\til{\alpha_i})$,
$\til{\alpha_i}^p = \til{a_i}$, and let $\til{D} = (\til{a},
\til{b})_\omega$. We will now demonstrate that $L = L_1 \otimes_K
\cdots \otimes_K L_r$, $\til{D}$, and $L \otimes_K \til{D}$ are
division algebras. To check this it suffices to show that these $K$
algebras are division algebras after being tensored over $K$ with
$\bQ_p(t)\(x_1\)\cdots\(x_{r-1}\)$. In other words, if we let $F =
\bQ_q(t)\(x_1\)\cdots\(x_{r-1}\)$, we may replace $L_i$ with
$L_i' = F(\til{\alpha})$ and $\til{D}$ with $D' = (\til{a}, \til{b})_{F,
\omega}$. Let $\cO_F = \bZ_p(t)[[x_1, \ldots, x_{r-1}]]$. Note that
$\cO_F$ is a valuation ring with quotient field $F$ and residue field
$\ov{F} = \bF_p(t) = \ov{K}$. If we let $\cO_{D'}$
be the algebra $(\til{a}, \til{b})_{\cO_F, \omega}$, and $\cO_{L_i'} =
\cO_F (\til{\alpha_i})$, then we have 
\begin{gather*}
\cO_{D'} \otimes_{\cO_F} F = D \ , \ \ \ \
\cO_{D'} \otimes_{\cO_F} \ov{F} = \cD \\
\cO_{L_i'} \otimes_{\cO_F} F = L_i' \ , \ \ \ \ 
\cO_{L_i'} \otimes_{\cO_F} \ov{F} = \ell_i \\
\big(\cO_{L_1'} \otimes_{\cO_F} \cdots \otimes_{\cO_F} \cO_{L_r'}\big)
\otimes_{\cO_F} F = L_1' \otimes_F \cdots \otimes_F L_r', \\
\big(\cO_{L_1'} \otimes_{\cO_F} \cdots \otimes_{\cO_F} \cO_{L_r'}\big)
\otimes_{\cO_F} \ov{F} = \ell_1 \otimes_{\ov{K}} \cdots
\otimes_{\ov{K}} \ell_r.
\end{gather*}

The fact that the desired algebras are division now follows from lemma
\ref{division_val}.

It remains to show that the division algebra $D = (\til{a}, \til{b})$
and the field $L = K(\til{\alpha_1}, \ldots, \til{\alpha_r})$ posses
unramified valuations. For each of the given algebras, this can be
checked directly. Let $\fB$ be the standard basis of one of these - in
the first case consisting of expressions of the form
$\alpha^i\beta^j$, $0 \leq i, j \leq p-1$ where $\alpha$ and $\beta$
are the standard generators of the symbol algebra, or in the second
case consisting of expressions of the form $\alpha_1^{n_1} \cdots
\alpha_r^{n_r}$, $0 \leq n_i \leq p-1$. We may then define our valuation via
$$v(\sum_{b \in \fB} \lambda_b b) = \underset{b \in \fB}{min}\{v(\lambda)\}.$$
One may now observe that this is an unramified valuation.
\end{proof}

\begin{lem} \label{division_val}
Let $R$ be a valuation ring with residue field $k$ and field of
fractions $F$. Suppose that $A$ is an $R$-algebra which is free of
finite rank as an $R$-module, and that $A \otimes_R k$ is a division
algebra. Then $A \otimes_R F$ is also a division algebra.
\end{lem}
\begin{proof}
The norm polynomial on $A$ defines a projective and hence proper
scheme $X$ over $Spec(R)$. It is easy to see that $A \otimes_R L$ is a
domain if and only if $X(L) = \emptyset$. In particular, if $A
\otimes_R F$ was not a domain, we would have a point in $X(F)$, and by
the valuative criterion for properness, we would thereby obtain a
point in $X(k)$, contradicting the fact that $A \otimes_R k$ is a
domian. The conclusion now follows from the observation that a finite
dimensional domain is a division algebra.
\end{proof}

\end{document}